\documentclass[12pt]{article}
\usepackage[latin1]{inputenc}
\usepackage[T1]{fontenc}
\usepackage{indentfirst}
\usepackage{amsmath}
\usepackage{amsfonts}
\usepackage[english]{babel}
\everymath{\displaystyle}

\font\ninerm=cmr10

\title{\Huge{A kinetic model for coagulation-fragmentation}}
\author{Damien BROIZAT}
\date{\small{Laboratoire J.A. Dieudonn\'e\\
Universit\'e de Nice-Sophia Antipolis, Parc Valrose,\\F-06108 Nice cedex 02, France.\\
E-mail : broizat@unice.fr}}

\begin{document}
\maketitle

\font\bba=msbm10 \font\bbb=msbm8
%\font\bbb=msbm10 scaled 800
\font\bbc=msbm6
%\font\bbc=msbm10 scaled 600
\newfam\bbfam
\textfont\bbfam=\bba \scriptfont\bbfam=\bbb
\scriptscriptfont\bbfam=\bbc
\def\bb{\fam\bbfam\bba}

\def\eps{\varepsilon}
\def\N{{\bb N}}
\def\Z{{\bb Z}}
\def\r{{\bb R}}
\def\C{{\bb C}}
\def\1{{1\hspace{-0.9mm}{\rm l}}}%\mbox{\normalsize I}}}
\def\emptyset{\hbox{$\displaystyle/\kern -5.97pt\circ$}}
\def\diam{\mbox{{\normalsize diam}}}
\def\supp{\mbox{{\normalsize supp}}}
\def\lip{\mbox{{\normalsize Lip}}}
\def\stackunder#1#2{\mathop{#1}\limits_{#2}}
\def\stackover#1#2{\mathop{#1}\limits^{#2}}
\def\CQFD{\unskip\kern 6pt\penalty 500%
\raise -1pt\hbox{\vrule\vbox to 8pt{\hrule width
6pt\vfill\hrule}\vrule}}
\def\ess{\mathop{\mbox{\rm ess}}}
\def\limsup{\mathop{\overline{\lim}}\limits}
\def\liminf{\mathop{\underline{\lim}}\limits}
\def\dv{\mathop{\mbox{\rm div}}\nolimits}
\def\sgn{\mathop{\mbox{\rm sgn}}}
\def\iint{\mathop{\int\mkern -12mu\int}}
\def\iiint{\mathop{\int\mkern -12mu\int\mkern -12mu\int}}
\def\iiiint{\mathop{\int\mkern -12mu\int\mkern -12mu\int\mkern -12mu\int}}
\def\tq{\mathrel{;\,}}
\def\pref#1{(\ref{#1})}
\def\pr{\noindent {\bf Proof. }}
\let\dsp=\displaystyle
\def\convol{\mathop{*}\limits}
\def\tr{\mathop{\mbox{\rm tr}}}
\def\meas{\mathop{\rm meas}}
\newcommand{\integr}[1]{\int_\r \int_\r #1  \,dx \,d\xi}
\newcommand{\integ}[1]{\int_\r  #1 \,d\xi}
\def\limd{\mathop{{\lim}}\limits}
\def\tod{\mathop{{\to}}\limits}
\def\supd{\mathop{{\sup}}\limits}
\def\maxd{\mathop{{\max}}\limits}
\def\mind{\mathop{{\min}}\limits}
\def\simd{\mathop{{\sim}}\limits}
\def\tw{\mathop{{\rightharpoonup}}\limits}

\newtheorem{Th}{Theorem}[section]
\newtheorem{Prop}[Th]{Proposition}
\newtheorem{Lemma}[Th]{Lemma}
\newtheorem{Defin}[Th]{Definition}
\newtheorem{Cor}[Th]{Corollary}
\newtheorem{Rk}{Remark}[section]
\newtheorem{Ex}{\sl Example}[section]

\def\theequation{\thesection.\arabic{equation}}
\def\thesection{\arabic{section}}
\def\thesubsection{\arabic{section}.\arabic{subsection}}
\def\thesubsubsection{\arabic{section}.\arabic{subsection}.\arabic{subsubsection}}
\newcommand\Section{%
\def\thesubsection{\arabic{section}}
\setcounter{Th}{0} \setcounter{Rk}{0} \setcounter{Ex}{0}
\setcounter{equation}{0}\section}
\newcommand\Subsection{%
\def\thesubsection{\arabic{section}.\arabic{subsection}}
\subsection}
\newcommand\Subsubsection{\subsubsection}

\baselineskip=12pt
\parskip=0pt plus 1pt
%\parskip=9.6pt plus 2pt minus 2pt

%\clearpage
%\tableofcontents
%\clearpage

\begin{center}
\textbf{Abstract}
\end{center}

The aim of this paper is to show an existence theorem for a kinetic model of coagulation-fragmentation with initial data satisfying the natural physical bounds, and assumptions of finite number of particles and finite $L^p$-norm. We use the notion of renormalized solutions introduced dy DiPerna and Lions in \cite{diperna-lions}, because of the lack of \textit{a priori} estimates. The proof is based on weak-compactness methods in $L^1$, allowed by $L^p$-norms propagation.

\begin{section}{Introduction}

Coalescence and fragmentation are general phenomena which appear in dynamics of particles, in various fields (polymers chemistry, raindrops formation, aerosols, ...). We can describe them at different scales, which lead to different mathematical points of view. First, we can study the dynamics at the microscopic level, with a system of $N$ particles which undergo successives mergers/break ups in a random way. We refer to the survey \cite{aldous} for this stochastic approach. Another way to describe coalescence and fragmentation is to consider the statistical properties of the system, introducing the statistical distribution of particles $f(t,m)$ of mass $m>0$ at time $t\geq0$ and studying its evolution in time. This approach is rather macroscopic. But we can put in an intermediate level, by considering a density $f$ which depends on more variables, like position $x$ or velocity $v$ of particles, and this description is more precise. Here, we start by discussing models with density, from the original (with $f=f(t,m)$) to the kinetic one (with $f=f(t,x,m,v)$), which is the setting of this work. 

Depending on the physical context, the mass variable is discrete (polymers formation) or continuous (raindrops formation). It leads to two sorts of mathematical models, with $m\in\mathbb{N}^{\star}$ or $m\in]0,+\infty[$, but we focus on the continuous case. To understand the relationship between discrete and continuous equations, see \cite{mischler4}.

\begin{subsection}{The original model}

The discrete equations of coagulation have been originally derived by Smoluchowski in \cite{smol1, smol2}, by studying the Brownian motion of colloidal particles. It had been extended to the continuous setting by M\"uller \cite{muller}, giving the following mathematical model, called the \textit{Smoluchowski's equation of coagulation}:

\begin{equation}\label{equation Smoluchowski}
\frac{\partial f}{\partial t}(t,m) =Q_c^+(f,f)-Q_c^-(f,f),\qquad (t,m)\in]0,+\infty[^2.
\end{equation}

\noindent This equation describes the evolution of the statistical mass distribution in time. At each time $t>0$, the term $Q_c^+(f,f)$ represents the gain of particles of mass $m$ created by coalescence between smaller ones, by the reaction $$\{m^{\star}\}+\{m-m^{\star}\}  \rightarrow  \{m\}.$$

\noindent The term $Q_c^-(f,f)$ is the depletion of particles of mass $m$ because of coagulation with other ones, following the reaction $$\{m\}+\{m^{\star}\}  \rightarrow  \{m+m^{\star}\}.$$ 

\noindent Namely, we have 

$$\left\{
\begin{array}{l}
Q_c^+(f,f)(t,m)=\frac{1}{2}\int_0^{m} A(m^{\star},m-m^{\star})f(t,m^{\star})f(t,m-m^{\star})dm^{\star},\\\\
Q_c^-(f,f)(t,m)=\int_0^{+\infty} A(m,m^{\star})f(t,m)f(t,m^{\star})dm^{\star},
\end{array}
\right.$$

\noindent where $A(m,m^{\star})$ is the coefficient of coagulation between two particles, which governs the frequency of coagulations, according to the mass of clusters. In his original model, Smoluchowski derived the following expression for $A$:

\begin{equation}\label{noyau originel}
A(m,m^{\star})=\left(m^{1/3}+{m^{\star}}^{1/3}\right)\left(m^{-1/3}+{m^{\star}}^{-1/3}\right).
\end{equation}
 
In many cases, coalescence is not the only mechanism governing the dynamics of particles, and other effects should be taken into account. A classical phenomenon which also occurs is the fragmentation of particles in two (or more) clusters, resulting from an internal dynamic (we do not deal here with fragmentation processes induced by particles collisions). This binary fragmentation is modeled by linear additional reaction terms in equation $(\ref{equation Smoluchowski})$, namely

$$\left\{
\begin{array}{l}
Q_f^+(f)(t,m)=\int_m^{+\infty} B(m',m)f(t,m')dm',\\\\
Q_f^-(f)(t,m)=\frac{1}{2}f(t,m)B_1(m),\qquad\mbox{where}\qquad B_1(m')=\int_0^{m'}B(m',m)dm.
\end{array}
\right.$$

\noindent The function $B(m',m)$ is the fragmentation kernel, it measures the frequency of the break-up of a mass $m'$ in two clusters $m$ and $m'-m$, for $m<m'$. So, at each time $t$, the term $Q_f^+(f)$ is the gain of particles of mass $m$, resulting from the following reaction of fragmentation: 
$$\{m'\}\rightarrow  \{m\}+\{m'-m\},$$

\noindent whereas $Q_f^-(f)$ stands for the loss of particles of mass $m$, because of a break-up into two smaller pieces, by the following way:
$$\{m\}\rightarrow  \{m^{\star}\}+\{m-m^{\star}\},\quad\mbox{with}\quad m^{\star}<m.$$

\noindent Thus, the continuous \textit{coagulation-fragmentation equation} writes

\begin{equation}\label{equation cofrag}
\frac{\partial f}{\partial t}(t,m) =Q_c^+(f,f)-Q_c^-(f,f)+Q_f^+(f)-Q_f^-(f),\qquad (t,m)\in]0,+\infty[^2.
\end{equation}

\noindent In the 90's, many existence and uniqueness results have been proved about this problem, see for instance \cite{stewart1, stewart2}, or \cite{lamb} for an approach by the semigroups of operators theory. These results are true under various growth hypotheses on kernels $A$ and $B$, but these assumptions often allow unbounded kernels, which is important from a physical point of view.

However, this coagulation-fragmentation model do not take the spatial distribution of particles into account. This leads to ``spatially inhomogeneous'' mathematical models, where the density of particles $f(t,x,m)$ depends also of a space variable $x\in\mathbb{R}^3$. 

\end{subsection}

\begin{subsection}{Spatially inhomogeneous models}

A first example consists of diffusive models, corresponding to the situation where particles follow a Brownian motion at the microscopic scale, with a positive and mass-dependent coefficient of diffusion $d(m)$. From a physical point of view, it implies that particles are sufficiently small to undergo the interaction with the medium, i.e the shocks with the molecules of the fluid in which the particles evolve. In the statistical description, a spatial-laplacian appears, giving the \textit{diffusive coagulation-fragmentation equation}:

\begin{equation}\label{equation cofrag diffusive}
\begin{array}{lll}
\frac{\partial f}{\partial t}(t,x,m)-d(m)\Delta_x f(t,x,m)=Q_c^+(f,f)-Q_c^-(f,f)+Q_f^+(f)-Q_f^-(f),\\\\
(t,x,m)\in]0,+\infty[\times\mathbb{R}^3\times]0,+\infty[.
\end{array}
\end{equation}
   
\noindent We refer to \cite{mischler3} for a global existence theorem for the discrete diffusive coagulation~-fragmentation equation in $L^1$, and to \cite{mischler2} for the continuous one, improved in \cite{mrr} (with less restrictive conditions on the kernels), then in \cite{amann} (with uniqueness of the solution).\vspace{5mm}

The second way to correct the spatially homogeneous problem is to assume that the particles are transported with a deterministic velocity $v$. At the statistical level, this adds a linear transport term $v.\nabla_x f$ to the equation $(\ref{equation cofrag})$. This velocity can be a given velocity $v=v(t,x,m)$ or the inner velocity of the particles. The first case has been studied in \cite{dubovski}, with an existence and uniqueness theorem, and furthermore the continuous dependance on the initial data. Physically, it corresponds to the dynamics of particles with rather low mass which follow a velocity drift depending only on the surrounding fluid. In the second case, particles are also identified by their momentum $p\in\mathbb{R}^3$ in addition to their mass $m$ (with $v=p/m$): we have a kinetic model, which is relevant to describe the dynamics of particles of varying size/mass according to coagulation/fragmentation events, like in aerosols. At the microscopic scale, the coagulation/fragmentation processes become ``multi-dimensional'', with mass-momentum conservation at each merger/break up according to the following scheme:

$$\begin{array}{lccc}
Coagulation: & \{m\}+\{m^{\star}\} & \rightarrow & \{m'\}\\
& \{p\}+\{p^{\star}\} & \rightarrow & \{p'\}\vspace{2mm}\\
Fragmentation: & \{m'\} & \rightarrow & \{m\}+\{m^{\star}\}\\
& \{p'\} & \rightarrow & \{p\}+\{p^{\star}\}\\
\end{array}$$

\noindent where $m':=m+m^{\star}$, $m>0$, $m^{\star}>0$, and $p':=p+p^{\star}.\vspace{2mm}$

\noindent Thus, in the statistical description, the density depends on time, position, mass and momentum: $f=f(t,x,m,p)$. But even if this kind of kinetic models provides a rather good description of phenomena, it is harder to study, so there are less results than for the diffusive ones. Moreover, it is difficult to know the exact physical form of the kernels. And finally, the numerical aspects are a real problem on these models: because of a high dimension (at least $7$ plus time), it seems to be very difficult, maybe impossible, to compute the solutions on a long time.

Concerning the results, a global existence theorem for the sole coagulation has been demonstrated in \cite{mischler1}. The proof is based on $L^p$-norms dissipation for any formal solution, and on weak-compactness methods in $L^1$. This result has been extended to a more general class of coalescence operators in \cite{jabin2} (but under stronger restriction on the initial data), with a very different method of proof. For the sole fragmentation, a difficulty is due to the blow-up of kinetic energy, which grows at each microscopic event. Thus, it is reasonable to take the internal energy of particles into account, which balances the gain of kinetic energy during a break up. With that modeling, the work \cite{jabin} provides global existence for a kinetic fragmentation model, with general growth assumptions on the kernel $B$, by using correct entropies.

The aim of this work is to combine both of these analysis. The main difficulty is the lack of \textit{a priori} estimates, because of the compensation problems between the two kernels: we are not able to define well the reaction term of coagulation with only the \textit{a priori} bounds (specifically, we can not say that $Q_c^-(f,f)$ lies in $L^1_{loc}$). That is why we use the DiPerna-Lions theory of renormalized solutions, introduced in \cite{diperna-lions} to show global existence for Boltzmann equation, which presents the same problems.

\end{subsection}

\begin{subsection}{Description of the kinetic model and outline of the paper}

Now, let us describe precisely the model we study. The parameters which describe the state of a particle are denoted by $$y:=(m,p,e)\in Y:=]0,+\infty[\times\mathbb{R}^3\times]0,+\infty[,$$ 

\noindent $m$ for the mass, $p$ the impulsion, and $e$ the internal energy. At the microscopic scale, coalescence and fragmentation conserve total energy  (kinetic energy + internal energy), thus we can compute the internal energy of daughter(s) particle(s).

$$\begin{array}{lccc}
Coagulation: & \{e\}+\{e^{\star}\} & \rightarrow & \{e'\}\vspace{2mm}\\
\end{array}$$ 

\noindent We have $$\frac{|p|^2}{2m}+e+\frac{|p^{\star}|^2}{2m^{\star}}+e^{\star}=\frac{|p+p^{\star}|^2}{2(m+m^{\star})}+e',$$ thus $$e'=e+e^{\star}+E_-(m,m^{\star},p,p^{\star}), \quad\mbox{where}\quad E_-(m,m^{\star},p,p^{\star}):=\frac{|m^{\star}p-mp^{\star}|^2}{2mm^{\star}(m+m^{\star})}\geq 0$$ 

\noindent ($E_-$ is the loss of kinetic energy resulting from the merger).

$$\begin{array}{lccc}
Fragmentation: & \{e'\} & \rightarrow & \{e\}+\{e^{\star}\}\vspace{2mm}\\
\end{array}$$ 

\noindent We have $$\frac{|p'|^2}{2m'}+e'=\frac{|p|^2}{2m}+e+\frac{|p'-p|^2}{2(m'-m)}+e^{\star},$$ thus 
$$e^{\star}=e'-e-E_+(m',m,p',p),\quad\mbox{where}\quad E_+(m',m,p',p):=\frac{|m'p-mp'|^2}{2mm'(m'-m)}\geq 0$$

\noindent ($E_+$ is the gain of kinetic energy resulting from the break up).\vspace{2mm}

\begin{Rk} Let us point out the following symmetries: 
$$E_-(m^{\star},m,p^{\star},p)=E_-(m,m^{\star},p,p^{\star})\quad\mbox{and}\quad E_+(m',m,p',p)=E_+(m',m'-m,p',p'-p),$$
and the relation: $E_-(m,m^{\star},p,p^{\star})=E_+(m+m^{\star},m,p+p^{\star},p)$, which is consistent with the two phenomena's reciprocity.\vspace{5mm}
\end{Rk}

\noindent We use the following notations:

\begin{itemize}
\item if $y=(m,p,e)$, $y^{\star}=(m^{\star},p^{\star},e^{\star})$, then we denote
$$y':=y+y^{\star}:=(m+m^{\star},p+p^{\star},e+e^{\star}+E_-(m,m^{\star},p,p^{\star})),$$ 
\item if $y=(m,p,e)$, $y'=(m',p',e')$, with $m<m'$ and $e<e'-E_+(m',m,p',p)$, then we say that
$y<y'$ and we denote
$$y^{\star}:=y'-y:=(m'-m,p'-p,e'-e-E_+(m',m,p',p)).$$
\end{itemize}

\noindent With this formalism, we naturally have $(y'-y)+y=y'$, but note carefully that $y<y'$ is not an order relation on $Y$.\vspace{5mm}

\begin{Rk} For all $y'\in Y,\quad\{y\in Y,\quad y<y'\}\subset]0,m'[\times B_{\sqrt{2m'e'+|p'|^2}}\times]0,e'[.$\\Denoting $Y_R:=]0,R[\times B_R\times]0,R[ \subset Y$, we have  
\begin{equation}\label{limitations}
y<y',\quad y'\in  Y_R\quad\Longrightarrow\quad y\in]0,R[\times B_{\sqrt{3}R}\times ]0,R[\subset Y_{2R}. 
\end{equation}
\end{Rk}\vspace{2mm}

\noindent Finally, we point out that the map $(m',m,p',p,e',e)\mapsto (m',m^{\star},p',p^{\star},e',e^{\star})$ is a diffeomorphism with $C^{\infty}$ regularity whithin the domain 
$$\{0<m<m',\,\, p,p'\in \mathbb{R}^3,\,\, 0<e<e'-E_+(m',m,p',p)\}\subset Y^2$$
which preserves volume.\vspace{5mm}

\noindent We denote by $f(t,x,m,p,e)=f(t,x,y)$ the particles density, which is a nonnegative function depending on time $t\geq0$, position $x\in\mathbb{R}^3$, and the mass-momentum-energy variable $y$. To shorten the notations, we set for each $t$,$x$, $f(y)=f(t,x,y)$, or ${f=f(t,x,y)}$, $f^{\star}=f(t,x,y^{\star})$, and $f'=f(t,x,y').$ The complete model then reads:

$$\left\{
\begin{array}{l}
\partial_t f +\frac{p}{m}.\nabla_x f=Q_c^+(f,f)-Q_c^-(f,f)+Q_f^+(f)-Q_f^-(f),\qquad (ECF)\\\\
t\in]0,+\infty[,\quad x\in\mathbb{R}^3,\quad y=(m,p,e)\in Y,
\end{array}
\right.$$

\noindent with\vspace{-4mm}

$$\hspace{-13mm}\left\{\begin{array}{l}
Q_c^+(f,f)(y)=\frac{1}{2}\int_Y A(y^{\star},y-y^{\star})f(y^{\star})f(y-y^{\star})\1_{\{y^{\star}< y\}}dy^{\star},\\\\
Q_c^-(f,f)(y)=f(y)Lf(y),\qquad Lf(y):=\int_Y A(y,y^{\star})f(y^{\star})dy^{\star},
\end{array}
\right.$$

\noindent and\vspace{-4mm}

$$\hspace{-13mm}\left\{\begin{array}{l}
Q_f^+(f)(y)=\int_Y B(y',y)f(y')\1_{\{y'> y\}}dy',\\\\
Q_f^-(f)(y)=\frac{1}{2}B_1(y)f(y),\qquad B_1(y'):=\int_Y B(y',y)\1_{\{y<y'\}}dy.
\end{array}
\right.$$

\noindent Functions $A$ et $B$ are respectively the coagulation and fragmentation kernels. They are nonnegative functions, independent of $(t,x)$, which satisfy the natural properties of symmetry:

\begin{equation}\label{symetrie A}
\forall (y,y^{\star})\in Y^2,\quad A(y,y^{\star})=A(y^{\star},y),
\end{equation}

\begin{equation}\label{symetrie B}
\forall (y,y')\in Y^2,\quad y<y',\quad B(y',y)=B(y',y^{\star}).
\end{equation}

\noindent The kernel $A(y,y^{\star})$ represents the coalescence rate between two particles $y$ and $y^{\star}$, whereas $B(y',y)$ is the fragmentation rate for a particle $y'$ which breaks in two clusters $y$ and $y^{\star}.$\vspace{5mm} 

\noindent We assume that $A$ fulfils the following structure assumption:

\begin{equation}\label{structure A}
\forall (y,y^{\star})\in Y^2,\quad A(y,y^{\star})\leq A(y,y')+A(y^{\star},y').
\end{equation}

\begin{Rk} We can insist on the fact that this assumption is more general than the classical \textit{Galkin-Tupchiev monotonicity condition}:
\begin{equation}\label{galkin}
\forall y<y^{\star},\quad A(y,y^{\star}-y)\leq A(y,y^{\star}).
\end{equation}
\noindent In the ``monodimensional'' case, the Smoluchowski kernel given by $(\ref{noyau originel})$ do not satisfy $(\ref{galkin})$ but satisfies $(\ref{structure A})$, that's why the first existence result established in \cite{mischler2} under Galkin-Tupchiev condition was extended in \cite{mrr} to kernels which satisfy $(\ref{structure A})$ only.
\end{Rk}

\noindent We also require that $A$ and $B$ have a mild growth:

\begin{equation}\label{asymptotique2 de A}
\forall R>0,\qquad \int_{Y_R}\frac{A(y,y^{\star})}{|y^{\star}|}dy\underset{|y^{\star}|\rightarrow +\infty}{\longrightarrow} 0,
\end{equation}

\begin{equation}\label{asymptotique de B}
\forall R>0,\qquad \int_{Y_R}\frac{B(y',y)}{|y'|}\1_{\{y<y'\}}dy\underset{|y'|\rightarrow +\infty}{\longrightarrow} 0,
\end{equation}

\noindent and $B$ is truncated as:

\begin{equation}\label{troncature du noyau de fragmentation}
\exists C_0>1,\quad 
\left\{
\begin{array}{lll}
& m'>C_0m\vspace{2mm}\\
\mbox{or} &\\ 
& e'+\frac{|p'|^2}{2m'}>C_0\left(e+\frac{|p|^2}{2m}\right)\\ 
\end{array}
\right.
\Longrightarrow \quad B(y',y)=0.
\end{equation}

\begin{Rk} The physical interpretation of this truncature assumption is to prevent the creation of too small clusters compared to the mother particle. From a mathematical point of view, it allows the total number of particles (the $L^1$-norm of $f$) to be finite at each time $t>0$. 
 
\end{Rk}

\noindent We also need to have $B_1$ locally bounded:

\begin{equation}\label{B_1 loc born\'ee}
\forall R>0,\qquad B_1\in L^{\infty}(Y_R),
\end{equation}

\noindent as well as $A$:

\begin{equation}\label{A loc born\'ee}
\forall R>0,\qquad A\in L^{\infty}(Y^2_R).
\end{equation}

\begin{Rk} Unfortunately, these assumptions of growth and boundedness are more restrictive, and in the monodimensional case, the Smoluchowski kernel $(\ref{noyau originel})$ doesn't satisfy them any more. The examples given in \cite{mischler1} for the sole coagulation, namely 
$$A(m,m^{\star},p,p^{\star})=(m^{\alpha}+{m^{\star}}^{\alpha})^2\left|\frac{p}{m}-\frac{p^{\star}}{m^{\star}}\right|,\quad 0\leq\alpha<1/2,$$ 

\noindent (for the dynamics of liquid droplets carried by a gaseous phase) or
 $$A(m,m^{\star},p,p^{\star})=\left(\frac{m+{m^{\star}}}{mm^{\star}}\right)^{\alpha}\left|\frac{p}{m}-\frac{p^{\star}}{m^{\star}}\right|^{\gamma},\quad 0\leq\alpha\leq 1,\quad -3<\gamma\leq 0,$$

\noindent (for a stellar dynamics context) do not fit neither. Here, we need coalescence kernels which are bounded when $m,m^{\star}\rightarrow 0$. But it is difficult to know the exact physical form of the kernels $A$ and $B$ because of the complexity of this kinetic model. Nevertheless, simple kernels given by $A(m,m^{\star})=m^{\alpha}+{m^{\star}}^{\alpha}$ with $0<\alpha<1$ fit. 
\end{Rk}

\noindent Finally, we assume that $A$ controls $B$ in the following sense:

$\exists s>1,\quad\exists 0<\delta<\frac{1}{6s-5}<1,$

\begin{equation}\label{comparaison des noyaux}
\forall y'\in Y,\quad \int_Y \frac{B(y',y)^s}{A(y,y')^{s-1}}\1_{\{y<y'\}}dy\leq 1+m'+\frac{|p'|^2}{2m'}+e'+\frac{1}{2}B_1(y')^{\delta}.
\end{equation}

\begin{Rk} This last assumption is more technical, but seems necessary to balance the contributions of the interaction terms $Q_c(f,f)$ and $Q_f(f)$, which are difficult to compare because $Q_c(f,f)$ is quadratic whereas $Q_f(f)$ is linear.
 
\end{Rk}

\noindent The paper consists in the proof of the following theorem.

\begin{Th}\label{Theoreme d'existence} Let $A$ and $B$ be kernels satisfying $(\ref{symetrie A})-(\ref{structure A})$ and $(\ref{asymptotique2 de A})-(\ref{comparaison des noyaux})$ and let $f^0$ be a nonnegative initial data which satisfies 

\begin{equation}\label{hypotheses sur la donn\'ee initiale}
K(f^0):=\int\!\!\!\int_{\mathbb{R}^3\times Y}\left(\left(1+m+\frac{|p|^2}{2m}+e+m|x|^2\right)f^0(x,y)+f^0(x,y)^s\right)dxdy<\infty,
\end{equation}

\noindent then for all $T>0$, there exists $f\in C([0,T],L^1(\mathbb{R}^3\times Y))$ such that $f(0)=f^0$ and $f$ is a renormalized solution to $(ECF)$. Moreover,

\begin{equation}\label{la solution est d'energie finie}
a.e\,\,t\in]0,T[,\qquad\int\!\!\!\int_{\mathbb{R}^3\times Y}\left(1+m+\frac{|p|^2}{2m}+e+m|x|^2\right)f(t,x,y)dxdy\leq K_T, 
\end{equation}

\begin{equation}\label{la solution est L^s}
a.e\,\,t\in]0,T[,\qquad\int\!\!\!\int_{\mathbb{R}^3\times Y}f(t,x,y)^sdxdy\leq K_T,  
\end{equation}

\noindent where the constant $K_T$ depends only on $C_0$, $T$, $K(f^0)$, $s$ and $\delta$ (defined in $(\ref{troncature du noyau de fragmentation})$ and $(\ref{comparaison des noyaux})$).
\end{Th}

\noindent Beyond existence problems, there are lots of others interesting subjects to explore. A first one concerns the mass conservation of the solution $f$, which is still an open problem for such kinetic models, even for the case of the sole coagulation. In the spatially homogeneous case, it has been shown in \cite{dubovski-stewart1} that total mass is preserved in time under mild growth hypotheses on kernels. But we know that in case of strong coagulation (typically the case of multiplicative kernels), a phenomenon of gelation occurs, which force the total mass of the system to decay from a certain time $T_g<+\infty$.
Then, problems of convergence to an equilibrium have been already studied for the spatially homogeneous equation $\cite{mischler5}$, under a detailed balance condition between kernels $A$ and $B$. We can also mention existence of self-similar solutions \cite{emrr, fournier, laurencot}, always for the spatially homogeneous case.\vspace{5mm}

In a first section, we will derive the \textit{a priori} estimates from the equation, giving the proper setting of the problem. Then, the proof of theorem is based on a well-known stability principle which says that if we are able to pass to the limit in the equation (the set of solutions is closed in a certain sense), then it would be easy to show the existence of a solution, applying the stability result to a sequence of approached problems which we can solve. So, the aim of the last section is to prove rigorously such a stability result and in fact, we work in the context of renormalized solutions, because the reaction term can not be defined as a distribution simply using the \textit{a priori} estimates. 

\end{subsection}

\begin{subsection}{Different notions of solutions}

\noindent We discuss here on different notions of solutions, recalling the DiPerna-Lions results. We set $Q(f,f)=Q_c^+(f,f) -Q_c^-(f,f) + Q_f^+ (f) - Q_f^-(f).$ 

\begin{Defin} Let $f$ be a nonnegative function, such that $f\in L^1_{loc}(]0,+\infty[\times\mathbb{R}^3\times Y)$. We say that $f$ is a renormalized solution of (ECF) if 
$$\frac{Q_c^{\pm}(f,f)}{1+f}\in L^1_{loc}(]0,+\infty[\times\mathbb{R}^3\times Y),\qquad
\frac{Q_f^{\pm}(f)}{1+f}\in L^1_{loc}(]0,+\infty[\times\mathbb{R}^3\times Y),$$
\noindent and if the function $g:=\log(1+f)$ satisfies the renormalized equation
$$\partial_t g +\frac{p}{m}.\nabla_x g=\frac{Q(f,f)}{1+f}\qquad (ECFR)$$
\noindent in $\mathcal{D}'(]0,+\infty[\times\mathbb{R}^3\times Y)$.
\end{Defin} 

\noindent The renormalization makes passing to the limit impossible because of the quotients in the reaction term, that is why we also need another notion of solution: the mild solutions, which only require local integrability in time and provide Duhamel's integral formulations to the problem in which we are able to pass to the limit.

\begin{Defin} Let $f$ be a nonnegative function, such that $f\in L^1_{loc}(]0,+\infty[\times\mathbb{R}^3\times Y)$. We say that $f$ is a mild solution of (ECF) if for almost all $(x,y)\in\mathbb{R}^3\times Y$,

$$\forall T>0,\qquad Q_c^{\pm}(f,f)^{\sharp}(t,x,y)\in L^1(]0,T[),\qquad Q_f^{\pm}(f)^{\sharp}(t,x,y)\in L^1(]0,T[),$$

\noindent and
\begin{equation}\label{mild solution}
\forall 0<s<t<\infty,\quad f^{\sharp}(t,x,y)-f^{\sharp}(s,x,y)=\int_s^t Q(f,f)^{\sharp}(\sigma,x,y)d\sigma,
\end{equation}
\noindent where $h^{\sharp}$ denotes the restriction to the characteristic lines of the equation:
$$h^{\sharp}(t,x,m,p,e):=h(t,x+t\frac{p}{m},m,p,e).$$
\end{Defin}

\noindent The following results are proved in \cite{diperna-lions}:

\begin{Lemma}$ $\vspace{2mm}

 (i) If $Q_c^{\pm}(f,f)\in L^1_{loc}(]0,+\infty[\times\mathbb{R}^3\times Y)$ and $Q_f^{\pm}(f)\in L^1_{loc}(]0,+\infty[\times\mathbb{R}^3\times Y),$ then the following assertions are equivalent:
\begin{itemize}
\item $f$ is a solution of $(ECF)$ in the sense of distributions.
\item $f$ is a renormalized solution of $(ECF)$.
\item $f$ is a mild solution of $(ECF)$.
\end{itemize}\vspace{2mm}

 (ii) If $f$ is a renormalized solution of $(ECF)$, then for all function $\beta\in C^1([0,+\infty[)$ such that $|\beta'(u)|\leq \frac{C}{1+u}$, the composed function $\beta(f)$ is a solution of
$$\partial_t{\beta(f)}+\frac{p}{m}.\nabla_x{\beta(f)}=\beta'(f)Q(f,f).$$
\noindent in the sense of distributions (here, the right side lies in $L^1_{loc}(]0,+\infty[\times\mathbb{R}^3\times Y)$).\vspace{2mm}

(iii) $f$ is a renormalized solution of $(ECF)$ if and only if\\ $f$ is a mild solution of $(ECF)$, $\frac{Q_c^{\pm}(f,f)}{1+f}\in L^1_{loc}(]0,+\infty[\times\mathbb{R}^3\times Y)$ and
$\frac{Q_f^{\pm}(f)}{1+f}\in L^1_{loc}(]0,+\infty[\times\mathbb{R}^3\times Y).$

\end{Lemma}

\end{subsection}

\end{section}

\begin{section}{\textit{A priori} estimates}
\noindent We consider the Cauchy problem 

\begin{equation}\label{probl\`eme de Cauchy}
\left\{
\begin{array}{ll}
(ECF)\\
f(0,x,y)=f^0(x,y).\\
\end{array}\right.
\end{equation}

\noindent We suppose in this section that (\ref{probl\`eme de Cauchy}) admit a sufficiently smooth solution $f$ in order to handle some formal quantities which are conserved or propagated by the equation $(ECF)$. More precisely, we will show the propagation of $L^q$ bounds for the solution along time:
\begin{Prop} If the initial data $f^0$ satisfies 

\begin{equation}\label{bornes sur la donn\'ee initiale}
K(f^0):=\int\!\!\!\int_{\mathbb{R}^3\times Y}\left(\left(1+m+\frac{|p|^2}{2m}+e+m|x|^2\right)f^0(x,y)+f^0(x,y)^s\right)dxdy<\infty,
\end{equation}

\noindent then for all $T>0$, any classical solution of the Cauchy problem (\ref{probl\`eme de Cauchy}) satisfies

\begin{equation}\label{bornes L^1, L^s et moments}
\sup_{t\in[0,T]}\int\!\!\!\int_{\mathbb{R}^3\times Y}\left(\left(1+m+\frac{|p|^2}{2m}+e+m|x|^2\right)f(t,x,y)+f(t,x,y)^q\right)dxdy\leq K_T,
\end{equation}

\noindent for all the exponents $q\in\,]\,5/6\,,\,s\,],$ and also

\begin{equation}\label{autres quantit\'es born\'ees}
\int_0^T\!\!\!\int_{\mathbb{R}^3} \left(D_1(f(t,x))+D_2(f(t,x))\right)dxdt\leq K_T,
\end{equation}

\noindent where

\begin{equation}\label{def D_1}
D_1(f(t,x)):=\frac{1}{2}\int\!\!\!\int_{Y\times Y} A(y,y^{\star})\sup(f,f^{\star})\inf(f,f^{\star})^s dy^{\star}dy\geq 0, 
\end{equation}

\begin{equation}\label{def D_2}
D_2(f(t,x)):=\frac{s-\delta}{2}\int\!\!\!\int_{Y\times Y} B(y',y)f(t,x,y')^s \1_{\{y<y'\}}dy'dy\geq 0, 
\end{equation}

\noindent and the constant $K_T$ depends only on $C_0$, $T$, $K(f^0)$, $s$ and $\delta$.
\end{Prop}

\begin{subsection}{Basic physical estimates}
\noindent We start with a fundamental formula, which gives the variation in time of some integral quantities involving the solution $f$. 

\begin{Lemma}
Let $H(u)$ be a function with $C^1$ regularity on $[0,+\infty[$ and $\Phi(y)$ a real or vectorial function.
We have 

\begin{equation}\label{formule fondamentale}
\begin{array}{lll}
& & \frac{d}{dt}\int\!\!\!\int_{\mathbb{R}^3\times Y}\Phi(y)H(f(t,x,y))dxdy\\\\
= & & \frac{1}{2}\int\!\!\!\int\!\!\!\int_{\mathbb{R}^3\times Y\times Y} Aff^{\star}\left(\Phi' d_u H (f')-\Phi d_u H(f)-\Phi^{\star} d_u H(f^{\star})\right)dy^{\star}dydx\\\\
& + & \frac{1}{2}\int\!\!\!\int\!\!\!\int_{\mathbb{R}^3\times Y\times Y} B f'\left(\Phi d_u H (f)+\Phi^{\star}d_u H(f^{\star})-\Phi' d_u H(f')\right)\1_{\{y<y'\}}dydy'dx,\\
\end{array} 
\end{equation}

\noindent where $d_u H=\frac{dH}{du}.$
\end{Lemma}

\noindent\underline{Proof}: Using $(ECF)$, we have
$$\begin{array}{lllll}
\frac{d}{dt}\int\!\!\!\int_{\mathbb{R}^3\times Y}\Phi(y)H(f)dxdy & = & & \hspace{-5mm}
\int\!\!\!\int_{\mathbb{R}^3\times Y}\Phi(y)\,d_u H(f)\,\partial_t f dydx\\\\
& = & & \hspace{-5mm}\int\!\!\!\int_{\mathbb{R}^3\times Y}\Phi(y)\,d_u H(f)\,\left(Q_c^+(f,f)-Q_c^-(f,f)\right)dydx\\\\
& & + & \int\!\!\!\int_{\mathbb{R}^3\times Y}\Phi(y)\,d_u H(f)\,\left(Q_f^+(f)-Q_f^-(f)\right)dydx\\\\
& & - & \int\!\!\!\int_{\mathbb{R}^3\times Y}div_x(-\Phi(y)H(f)\frac{p}{m})dydx.
\end{array}$$

\noindent The integral with divergence vanishes thanks to Stokes' formula. Whence
$$\begin{array}{lllllll}
\frac{d}{dt}\int\!\!\!\int_{\mathbb{R}^3\times Y}\Phi(y)H(f)dxdy
 & = & &\hspace{-5mm} \int\!\!\!\int_{\mathbb{R}^3\times Y}\Phi(y)\,d_u H(f)\,\left(Q_c^+(f,f)-Q_c^-(f,f)\right)dydx\\\\
 & & + & \int\!\!\!\int_{\mathbb{R}^3\times Y}\Phi(y)\,d_u H(f)\,\left(Q_f^+(f)-Q_f^-(f)\right)dydx.
\end{array}$$

\noindent Using Fubini's theorem (formally), we can write

$$\begin{array}{lllllll}
& & \hspace{-15mm}\frac{d}{dt}\int\!\!\!\int_{\mathbb{R}^3\times Y}\Phi(y)H(f)dxdy\\\\
= & & \hspace{-5mm}\frac{1}{2}\int\!\!\!\int\!\!\!\int_{\mathbb{R}^3\times Y\times Y}\Phi(y)\,d_u H(f)\, A(y^{\star},y-y^{\star})f(y^{\star})f(y-y^{\star})\1_{\{y^{\star}< y\}}dy^{\star}dydx\\\\
& - & \int\!\!\!\int\!\!\!\int_{\mathbb{R}^3\times Y\times Y}\Phi(y)\,d_u H(f) A(y,y^{\star})f(y)f(y^{\star})dy^{\star}dydx\\\\
& + & \int\!\!\!\int\!\!\!\int_{\mathbb{R}^3\times Y\times Y}\Phi(y)\,d_u H(f)\, B(y',y)f(y')\1_{\{y'> y\}}dy'dydx\\\\
& - & \frac{1}{2}\int\!\!\!\int_{\mathbb{R}^3\times Y}\Phi(y')\,d_u H(f') B_1(y')f(y')dy'dx.
\end{array}$$

\noindent If we change variables $(y^{\star},y-y^{\star})\rightarrow (y^{\star},y)$ in the first integral, we obtain

$$\begin{array}{lllllll}
& & \hspace{-15mm}\frac{d}{dt}\int\!\!\!\int_{\mathbb{R}^3\times Y}\Phi(y)H(f)dxdy\\\\
= & & \hspace{-5mm}\frac{1}{2}\int\!\!\!\int\!\!\!\int_{\mathbb{R}^3\times Y\times Y}\Phi(y+y^{\star})\,d_u H(f(y+y^{\star}))\, A(y^{\star},y)f(y^{\star})f(y)dy^{\star}dydx\\\\
& - & \int\!\!\!\int\!\!\!\int_{\mathbb{R}^3\times Y\times Y}\Phi(y)\,d_u H(f) A(y,y^{\star})f(y)f(y^{\star})dy^{\star}dydx\\\\
& + & \int\!\!\!\int\!\!\!\int_{\mathbb{R}^3\times Y\times Y}\Phi(y)\,d_u H(f)\, B(y',y)f(y')\1_{\{y'> y\}}dy'dydx\\\\
& - & \frac{1}{2}\int\!\!\!\int\!\!\!\int_{\mathbb{R}^3\times Y\times Y}\Phi(y')\,d_u H(f')\, B(y',y)f(y')\1_{\{y<y'\}}dydy'dx.
\end{array}$$

\noindent The symmetry of $A$ allows us to write

$$\begin{array}{llll}
& & & \hspace{-15mm}\int\!\!\!\int\!\!\!\int_{\mathbb{R}^3\times Y\times Y}\Phi(y)\,d_u H(f) A(y,y^{\star})f(y)f(y^{\star})dy^{\star}dydx\\\\
& = & & \hspace{-5mm}\frac{1}{2}\int\!\!\!\int\!\!\!\int_{\mathbb{R}^3\times Y\times Y}\Phi(y)\,d_u H(f) A(y,y^{\star})f(y)f(y^{\star})dy^{\star}dydx\\\\
& & + &\frac{1}{2}\int\!\!\!\int\!\!\!\int_{\mathbb{R}^3\times Y\times Y}\Phi(y^{\star})\,d_u H(f^{\star})A(y,y^{\star})f(y)f(y^{\star})dy^{\star}dydx,
\end{array}$$

\noindent using the change of variables $(y,y^{\star})\rightarrow (y^{\star},y)$.

\noindent The same applies to $B$ with $(y',y)\rightarrow (y',y'-y)$:
$$\begin{array}{llll}
& & & \hspace{-15mm}\int\!\!\!\int\!\!\!\int_{\mathbb{R}^3\times Y\times Y}\Phi(y)\,d_u H(f)\, B(y',y)f(y')\1_{\{y'> y\}}dy'dydx\\\\
& = & &\hspace{-5mm}\frac{1}{2}\int\!\!\!\int\!\!\!\int_{\mathbb{R}^3\times Y\times Y}\Phi(y)\,d_u H(f)\, B(y',y)f(y')\1_{\{y'> y\}}dy'dydx\\\\
& & +& \frac{1}{2}\int\!\!\!\int\!\!\!\int_{\mathbb{R}^3\times Y\times Y}\Phi(y'-y)\,d_u H(f(y'-y))\, B(y',y)f(y')\1_{\{y'> y\}}dy'dydx.
\end{array}$$

\begin{flushright}
$\Box$ 
\end{flushright}\vspace{5mm}

\noindent Applying this lemma with $H(u)=u$, it gives
\begin{equation}\label{moments de f}
\begin{array}{lllll}
\frac{d}{dt}\int\!\!\!\int_{\mathbb{R}^3\times Y}\Phi(y)f dxdy & = & & \hspace{-5mm}\frac{1}{2}\int\!\!\!\int\!\!\!\int_{\mathbb{R}^3\times Y\times Y} Aff^{\star}\left(\Phi'-\Phi-\Phi^{\star}\right)dy^{\star}dydx\\\\
& & + & \frac{1}{2}\int\!\!\!\int\!\!\!\int_{\mathbb{R}^3\times Y\times Y} Bf'\left(\Phi+\Phi^{\star}-\Phi'\right)\1_{\{y<y'\}}dy'dydx.
\end{array}
\end{equation}

\noindent Choosing $\Phi(y)=m$, we obtain mass conservation:
\begin{equation}\label{conservation masse}
\frac{d}{dt}\int\!\!\!\int_{\mathbb{R}^3\times Y}mf(t,x,y)dxdy=0.
\end{equation}

\noindent With $\Phi(y)=p$, we also get the momentum conservation:
\begin{equation}\label{conservation moments}
\frac{d}{dt}\int\!\!\!\int_{\mathbb{R}^3\times Y}pf(t,x,y)dxdy=0.
\end{equation}

\noindent Then, choosing $\Phi(y)=\frac{|p|^2}{2m}+e$, we recover the total energy conservation:
\begin{equation}\label{conservation \'energie}
\frac{d}{dt}\int\!\!\!\int_{\mathbb{R}^3\times Y}\left(\frac{|p|^2}{2m}+e\right)f(t,x,y)dxdy=0.
\end{equation}

\noindent Moreover, we can control space momenta:

\begin{Lemma} For all $T>0$, there exists a constant $C_T>0$ that
\begin{equation}\label{moments en espace}
\forall t\in [0,T],\quad \int\!\!\!\int_{\mathbb{R}^3\times Y} m{|x|^2}f(t,x,y)dxdy\leq C_T.
\end{equation}
\end{Lemma}
 
\noindent\underline{Proof}: In view of the equation $(ECF)$ and the Stokes formula, we have
 
$$\begin{array}{lllll}
\frac{d}{dt}\int\!\!\!\int_{\mathbb{R}^3\times Y} m{|x|^2}fdxdy & = &  -\int\!\!\!\int_{\mathbb{R}^3\times Y} {|x|^2}p.\nabla_x fdxdy\\\\
& = & 2\int\!\!\!\int_{\mathbb{R}^3\times Y} x.p \,\,f(t,x,y)dxdy\\\\
& \leq & 2\left(\int\!\!\!\int_{\mathbb{R}^3\times Y} m{|x|^2}fdxdy\right)^{1/2}\left(\int\!\!\!\int_{\mathbb{R}^3\times Y} \frac{|p|^2}{m}fdxdy\right)^{1/2},
\end{array}$$
 
\noindent and we conclude with (\ref{conservation \'energie}) and Gronwall's lemma.
\begin{flushright}
$\Box$ 
\end{flushright}\vspace{5mm}

\noindent Finally, we can control the number of particles in finite time:

\begin{Lemma} \noindent We set
$$\begin{array}{lll}
N_0:=\int\!\!\!\int_{\mathbb{R}^3\times Y} f^0(x,y)dxdy,\quad M_0:=\int\!\!\!\int_{\mathbb{R}^3\times Y} mf^0(x,y)dxdy,\\\\
E_0:=\int\!\!\!\int_{\mathbb{R}^3\times Y}\left(\frac{|p|^2}{2m}+e\right)f^0(x,y)dxdy. &\\
\end{array}$$Then, there exists a constant $C>0$ depending only on $C_0$ that
\begin{equation}\label{bornes L^1 en temps fini}
\forall T>0,\quad\forall t\in [0,T],\quad \int\!\!\!\int_{\mathbb{R}^3\times Y} f(t,x,y)dxdy\leq (N_0+CT(M_0+E_0))e^{CT}+M_0+E_0.
\end{equation}
\end{Lemma}

\noindent\underline{Proof}: We use formula (\ref{moments de f}) with $\Phi(y)=\1_{\{m\leq 1,\,\, e+\frac{|p|^2}{2m}\leq 1\}}$. Since $\Phi$ is nonnegative and subadditive in the sense of coalescence (ie $\Phi'\leq\Phi+\Phi^{\star}$), we have

$$\begin{array}{llllllll}
\frac{d}{dt}\int\!\!\!\int_{\mathbb{R}^3\times Y}\1_{\{m\leq 1,\,\, e+\frac{|p|^2}{2m}\leq 1\}} fdydx
& \leq & \frac{1}{2}\int\!\!\!\int\!\!\!\int_{\mathbb{R}^3\times Y\times Y}\!\!Bf'\left(\Phi+\Phi^{\star}-\Phi'\right)\1_{\{y<y'\}}dydy'dx\\\\
& = & \int\!\!\!\int\!\!\!\int_{\mathbb{R}^3\times Y\times Y}Bf'\left(\Phi-\frac{\Phi'}{2}\right)\1_{\{y<y'\}}dydy'dx\\\\
& \leq & \int\!\!\!\int\!\!\!\int_{\mathbb{R}^3\times Y\times Y}Bf'\Phi\1_{\{y<y'\}}dydy'dx\\\\
& = & \int\!\!\!\int\!\!\!\int_{\mathbb{R}^3\times Y\times Y}Bf'\1_{\{y<y',\,\, m\leq 1,\,\, e+\frac{|p|^2}{2m}\leq 1\}}dydy'dx.\\\\
\end{array}$$

\noindent In the last integral, if $m'>C_0$, then, since $m\leq 1$, we have ${B(y',y)=0}$ according to assumption (\ref{troncature du noyau de fragmentation}). The same applies if $e'+\frac{|p'|^2}{2m'}>C_0$.
Thus,

$$\begin{array}{lllll}
& & & \hspace{-20mm}\frac{d}{dt}\int\!\!\!\int_{\mathbb{R}^3\times Y} f(t,x,y)\1_{\{m\leq 1,\,\, e+\frac{|p|^2}{2m}\leq 1\}}dydx\\\\
& \leq & & \int\!\!\!\int_{\mathbb{R}^3\times Y}\left(\int_{Y}B(y',y)\1_{\{y<y'\}}dy\right)f(t,x,y')\1_{\{m'\leq C_0,\,\,e'+\frac{|p'|^2}{2m'}\leq C_0\}}dy'dx\\\\
& = & & \int\!\!\!\int_{\mathbb{R}^3\times Y}B_1(y')f(t,x,y')\1_{\{m'\leq C_0,\,\,e'+\frac{|p'|^2}{2m'}\leq C_0\}}dy'dx.
\end{array}$$

\noindent Denoting 
$C:=\sup_{y'\in Y_{2C_0}}B_1(y')$, we obtain

$$\begin{array}{lllll}
& & & \hspace{-20mm}\frac{d}{dt}\int\!\!\!\int_{\mathbb{R}^3\times Y}f(t,x,y)\1_{\{m\leq 1,\,\, e+\frac{|p|^2}{2m}\leq 1\}}dydx\\\\
& \leq & & C\int\!\!\!\int_{\mathbb{R}^3\times Y}f(t,x,y)dydx\\\\
& \leq & & C\left(\int\!\!\!\int_{\mathbb{R}^3\times Y}f(t,x,y)\1_{\{m\leq 1,\,\, e+\frac{|p|^2}{2m}\leq 1\}}dydx+\int\!\!\!\int_{\mathbb{R}^3\times Y}mf(t,x,y)dydx\right. \\\\
& & & \hphantom{C}\left. +\int\!\!\!\int_{\mathbb{R}^3\times Y}\left(e+\frac{|p|^2}{2m}\right)f(t,x,y)dydx\right).\\\\
\end{array}$$

\noindent Using (\ref{conservation masse}) and (\ref{conservation \'energie}), we have 

$$\begin{array}{ll}
& \hspace{-20mm}\frac{d}{dt}\int\!\!\!\int_{\mathbb{R}^3\times Y}f(t,x,y)\1_{\{m\leq 1,\,\, e+\frac{|p|^2}{2m}\leq 1\}}dydx\\\\
\leq & C\int\!\!\!\int_{\mathbb{R}^3\times Y}f(t,x,y)\1_{\{m\leq 1,\,\, e+\frac{|p|^2}{2m}\leq 1\}}dydx+C(M_0+E_0).   
\end{array}$$

\noindent We integrate this inequality in time. Then, Gronwall's lemma provides  
$$\forall T>0,\quad\forall t\in [0,T],\quad \int\!\!\!\int_{\mathbb{R}^3\times Y}\!\!f(t,x,y)\1_{\{m\leq 1,\,\, e+\frac{|p|^2}{2m}\leq 1\}}dydx\leq (N_0+CT(M_0+E_0))e^{CT}.$$

\noindent We conclude noting that
$$\begin{array}{lllll}
\int\!\!\!\int_{\mathbb{R}^3\times Y}f(t,x,y)dydx & \leq & & \hspace{-5mm}\int\!\!\!\int_{\mathbb{R}^3\times Y}f(t,x,y)\1_{\{m\leq 1,\,\, e+\frac{|p|^2}{2m}\leq 1\}}dydx\\\\
& & + & \int\!\!\!\int_{\mathbb{R}^3\times Y}\left(m+e+\frac{|p|^2}{2m}\right)f(t,x,y)dydx,
\end{array}$$

\noindent and using (\ref{conservation masse}) and (\ref{conservation \'energie}) again.
\begin{flushright}
$\Box$ 
\end{flushright}

\noindent To summarize, if we set $E(x,y)=1+m+\frac{|p|^2}{2m}+e+m|x|^2$, and if we suppose that the initial data satisfies $$K(f^0):=\int\!\!\!\int_{\mathbb{R}^3\times Y} E(x,y)f^0(x,y)dxdy<+\infty,$$
\noindent then, for all $T>0$, there exists a constant $K_T$ (depending on $T$, $C_0$ and $K(f^0)$) such that
\begin{equation}\label{premier r\'esum\'e des estimations}
\sup_{t\in[0,T]}\int\!\!\!\int_{\mathbb{R}^3\times Y} E(x,y)f(t,x,y)dxdy\leq K_T. 
\end{equation}

\begin{Rk} For $\gamma>5$, we have
\begin{equation}\label{int\'egrabilit\'e de l'\'energie}
\int\!\!\!\int_{\mathbb{R}^3\times Y}\frac{1}{E^{\gamma}(x,y)}dxdy<+\infty.
\end{equation}
\noindent It will be very useful to show that some $L^{q}$ bounds of $f$ (for $5/6< q<1$ and $q=s>1$) also propagate in time.\vspace{2mm}
\end{Rk}

\end{subsection}

\begin{subsection}{$L^q$ bounds} 

\noindent Obtaining $L^q$ bounds propagation is more technical, that is why we split the proof in several lemmas.  

\begin{Lemma} Let $\beta\in\,]\,5/6\,,\,1\,[$. Then, for all $T>0$, there exists a constant $K_T$ (depending on $T$, $C_0$ and $K(f^0)$) such that
\begin{equation}\label{bornes L^q, q<1}
\sup_{t\in[0,T]}\int\!\!\!\int_{\mathbb{R}^3\times Y} f^{\beta}(t,x,y)dxdy\leq K_T. 
\end{equation}
\end{Lemma}

\noindent\underline{Proof}: Writing
$$\int\!\!\!\int_{\mathbb{R}^3\times Y} f^{\beta}(t,x,y)dxdy=\int\!\!\!\int_{\mathbb{R}^3\times Y} f^{\beta}(t,x,y)\frac{E^{\beta}(x,y)}{E^{\beta}(x,y)}dxdy,$$
\noindent we use Young inequality 
\begin{equation}\label{Young facile}
\forall \alpha>1,\quad\forall u\geq 0,\quad\forall v\geq 0,\qquad uv\leq\frac{u^{\alpha}}{\alpha}+\frac{v^{\alpha^{\star}}}{\alpha^{\star}}
\end{equation}

\noindent with $u=f^{\beta}(t,x,y)E^{\beta}(x,y)$, $v=\frac{1}{E^{\beta}(x,y)}$, and $\alpha=\frac{1}{\beta}>1,$

\noindent and obtain
$$\begin{array}{llll}
\int\!\!\!\int_{\mathbb{R}^3\times Y} f^{\beta}(t,x,y)dxdy & \leq & & \!\!\!\!\!\!\beta \int\!\!\!\int_{\mathbb{R}^3\times Y} E(x,y)f(t,x,y)dxdy\\\\
& & + & (1-\beta)\int\!\!\!\int_{\mathbb{R}^3\times Y} \frac{1}{E^{\frac{\beta}{1-\beta}}(x,y)}dxdy.   
\end{array}
$$

\noindent We conclude with (\ref{premier r\'esum\'e des estimations}) and (\ref{int\'egrabilit\'e de l'\'energie}).
\begin{flushright}
$\Box$ 
\end{flushright}

\begin{Lemma} For any convex and nonnegative function $H\in C^1([0,+\infty[)$ such that $H(0)=0$, and for all $t>0$, we have  
\begin{equation}\label{bornes en fonction du moment B2 f}
\begin{array}{lllll}
\int\!\!\!\int_{\mathbb{R}^3\times Y} H(f(t,x,y))dxdy  & \leq & & \hspace{-5mm}\int\!\!\!\int_{\mathbb{R}^3\times Y} H(f^0(x,y))dxdy\\\\ 
& & - & \frac{1}{2}\int_0^t\!\!\!\int\!\!\!\int\!\!\!\int_{\mathbb{R}^3\times Y\times Y} \!\!\!\!\!\!\!A\sup(f,f^{\star})H(\inf(f,f^{\star}))dy^{\star}dydxd\tau\\\\
&  & + & \int_0^t\!\!\!\int\!\!\!\int\!\!\!\int_{\mathbb{R}^3\times Y\times Y} A' H\left(\frac{B}{A'}\right) f'\1_{\{y<y'\}}dydy'dxd\tau\\\\
&  & - & \frac{1}{2}\int_0^t\!\!\!\int\!\!\!\int\!\!\!\int_{\mathbb{R}^3\times Y\times Y} B f'd_u H(f')\1_{\{y<y'\}}dydy'dxd\tau, 
\end{array}
\end{equation}

\noindent where $A=A(y,y^{\star})$, $A'=A(y,y')$, $B=B(y',y)$.
\end{Lemma}

\noindent\underline{Proof}: The formula (\ref{formule fondamentale}) yields

\begin{equation}\label{d\'erivation de H sans sym\'etrie}
\begin{array}{lllll}
\frac{d}{dt}\int\!\!\!\int_{\mathbb{R}^3\times Y} H(f)dxdy 
& = & & \hspace{-5mm}\int\!\!\!\int\!\!\!\int_{\mathbb{R}^3\times Y\times Y} Aff^{\star}\left(\frac{d_u H (f')}{2}-d_u H(f)\right)dy^{\star}dydx\\\\
& & + & \int\!\!\!\int\!\!\!\int_{\mathbb{R}^3\times Y\times Y} \!B f'\left(d_u H (f)-\frac{d_u H(f')}{2}\right)\1_{\{y<y'\}}dydy'dx.\\
\end{array}
\end{equation}

\noindent Let us rewrite the term $I_1:=\int\!\!\!\int\!\!\!\int_{\mathbb{R}^3\times Y\times Y}\! Aff^{\star}d_u H(f')dy^{\star}dydx,$ by the following way:
$$I_1=\int\!\!\!\int\!\!\!\int_{\mathbb{R}^3\times Y\times Y} A\inf(f,f^{\star})\sup(f,f^{\star})d_u H(f')dy^{\star}dydx.$$
We use the Young inequality:

\begin{equation}\label{inegalit\'e convexe conjugu\'ee}
\forall u>0,\quad \forall v>0,\quad uv\leq H(u)+H^\star(v)
\end{equation}
with $u=\sup(f,f^{\star})$ and $v=d_u H(f')$, where $H^{\star}$ stands for the convex conjugate function of $H$. A simple calculus shows that 
$$H^{\star}(d_u H(u))=ud_u H(u)-H(u),$$

\noindent and this quantity is nonnegative, by the assumptions on $H$. We denote $$\Theta(u):=H^{\star}(d_u H(u))\geq 0.$$

\noindent It leads to the inequality:
$$\begin{array}{llll}
I_1 & \leq & & \hspace{-5mm}\int\!\!\!\int\!\!\!\int_{\mathbb{R}^3\times Y\times Y} A\inf(f,f^{\star})H(\sup(f,f^{\star}))dy^{\star}dydx\\\\
& & + & \int\!\!\!\int\!\!\!\int_{\mathbb{R}^3\times Y\times Y}A\inf(f,f^{\star})\Theta(f')dy^{\star}dydx.
\end{array}$$

\noindent We can dominate the second term of the right member using the hypothesis (\ref{structure A}) by the following way:

$$\begin{array}{llll}
& & &\hspace{-12mm}\int\!\!\!\int\!\!\!\int_{\mathbb{R}^3\times Y\times Y}A\inf(f,f^{\star})\Theta(f')dy^{\star}dydx\\\\
 & \leq & & \hspace{-5mm}\int\!\!\!\int\!\!\!\int_{\mathbb{R}^3\times Y\times Y} (A(y,y+y^{\star})+A(y^{\star},y+y^{\star}))\inf(f,f^{\star})\Theta(f')dy^{\star}dydx\\\\
& \leq & & \hspace{-5mm}\int\!\!\!\int\!\!\!\int_{\mathbb{R}^3\times Y\times Y}\!\!\!\!\!A(y,y+y^{\star})f\Theta(f')dy^{\star}dydx+\int\!\!\!\int\!\!\!\int_{\mathbb{R}^3\times Y\times Y}\!\!\!\!\!A(y^{\star},y+y^{\star})f^{\star}\Theta(f')dy^{\star}dydx\\\\
& = & & \hspace{-5mm}2\int\!\!\!\int\!\!\!\int_{\mathbb{R}^3\times Y\times Y}A(y^{\star},y+y^{\star})f(y^{\star})\Theta(f(y+y^{\star}))dy^{\star}dydx\\\\
& = & & \hspace{-5mm}2\int\!\!\!\int\!\!\!\int_{\mathbb{R}^3\times Y\times Y}A(y^{\star},y)f(y^{\star})\Theta(f(y))\1_{\{y^{\star}< y\}}dy^{\star}dydx
\end{array}$$ 

\noindent (the last identity resulting from the change of variables $(y^{\star},y+y^{\star})\rightarrow(y^{\star},y))$.\vspace{2mm}

\noindent This yields
$$\begin{array}{llll}
I_1 & \leq & & \hspace{-5mm}\int\!\!\!\int\!\!\!\int_{\mathbb{R}^3\times Y\times Y} A\inf(f,f^{\star})H(\sup(f,f^{\star}))dy^{\star}dydx\\\\
& & + & 2\int\!\!\!\int\!\!\!\int_{\mathbb{R}^3\times Y\times Y}Af^{\star}\Theta(f)\1_{\{y^{\star}< y\}}dy^{\star}dydx.
\end{array}$$ 

\noindent Thus, we have the following control of the coagulation contribution in (\ref{d\'erivation de H sans sym\'etrie}):

$$\begin{array}{llll}
& & & \hspace{-20mm}\int\!\!\!\int\!\!\!\int_{\mathbb{R}^3\times Y\times Y} Aff^{\star}\left(\frac{d_u H (f')}{2}-d_u H(f)\right)dy^{\star}dydx \\\\
\leq & & & \hspace{-5mm}\frac{1}{2}\int\!\!\!\int\!\!\!\int_{\mathbb{R}^3\times Y\times Y} A\inf(f,f^{\star})H(\sup(f,f^{\star}))dy^{\star}dydx\\\\
& & + & \int\!\!\!\int\!\!\!\int_{\mathbb{R}^3\times Y\times Y}Af^{\star}\Theta(f)\1_{\{y^{\star}< y\}}dy^{\star}dydx
 - \int\!\!\!\int\!\!\!\int_{\mathbb{R}^3\times Y\times Y}Aff^{\star}d_uH(f)dy^{\star}dydx.
\end{array}$$

\noindent Now, we can write the right member of this inequality by the following way:

$$\begin{array}{llll}
= & &  - & \frac{1}{2}\int\!\!\!\int\!\!\!\int_{\mathbb{R}^3\times Y\times Y} A\sup(f,f^{\star})H(\inf(f,f^{\star}))dy^{\star}dydx\\\\
& & + &\frac{1}{2}\int\!\!\!\int\!\!\!\int_{\mathbb{R}^3\times Y\times Y}AfH(f^{\star})dy^{\star}dydx
+\frac{1}{2}\int\!\!\!\int\!\!\!\int_{\mathbb{R}^3\times Y\times Y}Af^{\star}H(f)dy^{\star}dydx\\\\
& & + & \int\!\!\!\int\!\!\!\int_{\mathbb{R}^3\times Y\times Y}Af^{\star}\Theta(f)\1_{\{y^{\star}< y\}}dy^{\star}dydx
 - \int\!\!\!\int\!\!\!\int_{\mathbb{R}^3\times Y\times Y}Aff^{\star}d_uH(f)dy^{\star}dydx\\\\
= & &  - & \frac{1}{2}\int\!\!\!\int\!\!\!\int_{\mathbb{R}^3\times Y\times Y} A\sup(f,f^{\star})H(\inf(f,f^{\star}))dy^{\star}dydx\\\\
& & + &\int\!\!\!\int\!\!\!\int_{\mathbb{R}^3\times Y\times Y}Af^{\star}\Theta(f)\1_{\{y^{\star}< y\}}dy^{\star}dydx-\int\!\!\!\int\!\!\!\int_{\mathbb{R}^3\times Y\times Y}Af^{\star}\Theta(f)dy^{\star}dydx.
\end{array}$$

\noindent We deduce

\begin{equation}\label{contribution coagulation}
\begin{array}{llll}
& & & \hspace{-20mm}\int\!\!\!\int\!\!\!\int_{\mathbb{R}^3\times Y\times Y} Aff^{\star}\left(\frac{d_u H (f')}{2}-d_u H(f)\right)dy^{\star}dydx \\\\
\leq & &  - & \frac{1}{2}\int\!\!\!\int\!\!\!\int_{\mathbb{R}^3\times Y\times Y} A\sup(f,f^{\star})H(\inf(f,f^{\star}))dy^{\star}dydx\\\\
& & + &\int\!\!\!\int\!\!\!\int_{\mathbb{R}^3\times Y\times Y}Af^{\star}\Theta(f)\1_{\{y^{\star}< y\}}dy^{\star}dydx-\int\!\!\!\int\!\!\!\int_{\mathbb{R}^3\times Y\times Y}Af^{\star}\Theta(f)dy^{\star}dydx. 
\end{array}
\end{equation}

\noindent Then, we can also control the fragmentation contribution:
$$\begin{array}{lll}
\int\!\!\!\int\!\!\!\int_{\mathbb{R}^3\times Y\times Y} B f'd_u H (f)\1_{\{y<y'\}}dydy'dx
-\int\!\!\!\int\!\!\!\int_{\mathbb{R}^3\times Y\times Y} B f'\frac{d_u H(f')}{2}\1_{\{y<y'\}}dydy'dx.\\
\end{array}$$

\noindent We rewrite the first term by the following way:

$$\begin{array}{lllll}
& \int\!\!\!\int\!\!\!\int_{\mathbb{R}^3\times Y\times Y} B f'd_u H (f)\1_{\{y<y'\}}dydy'dx 
=  \int\!\!\!\int\!\!\!\int_{\mathbb{R}^3\times Y\times Y} \frac{B}{A'} A' f'd_u H (f)\1_{\{y<y'\}}dydy'dx
\end{array}$$

\noindent and we use (\ref{inegalit\'e convexe conjugu\'ee}) again, with $u=\frac{B}{A'}$ and $v=d_u H(f)$, whence 

\begin{equation}\label{contribution fragmentation}
\begin{array}{lll}
& & \hspace{-15mm}\int\!\!\!\int\!\!\!\int_{\mathbb{R}^3\times Y\times Y} B f'\left(d_u H (f)-\frac{d_u H(f')}{2}\right)\1_{\{y< y'\}}dydy'dx\\\\
\leq & & \int\!\!\!\int\!\!\!\int_{\mathbb{R}^3\times Y\times Y}\!\! H\left(\frac{B}{A'}\right)A' f'\1_{\{y<y'\}}dydy'dx +  \int\!\!\!\int\!\!\!\int_{\mathbb{R}^3\times Y\times Y}\!\!\! A' f'\Theta(f)\1_{\{y<y'\}}dydy'dx\\\\
& - & \int\!\!\!\int\!\!\!\int_{\mathbb{R}^3\times Y\times Y} B f'\frac{d_u H(f')}{2}\1_{\{y<y'\}}dydy'dx.\\\\
\end{array}
\end{equation}

\noindent Eventually, using (\ref{d\'erivation de H sans sym\'etrie}), (\ref{contribution coagulation}) and (\ref{contribution fragmentation}), we infer

\begin{equation}\label{estimation sur la d\'eriv\'ee}
\begin{array}{llll}
\frac{d}{dt}\int\!\!\!\int_{\mathbb{R}^3\times Y} H(f(t,x,y))dxdy & \leq & - & \frac{1}{2}\int\!\!\!\int\!\!\!\int_{\mathbb{R}^3\times Y\times Y} \!\!A\sup(f,f^{\star})H(\inf(f,f^{\star}))dy^{\star}dydx\\\\
& & - & \int\!\!\!\int\!\!\!\int_{\mathbb{R}^3\times \left(Y^2-\left(\{y< y^{\star}\}\cup\{y^{\star}< y\}\right)\right)}\! A f^{\star}\Theta(f)dydy^{\star}dx\\\\
& & + & \int\!\!\!\int\!\!\!\int_{\mathbb{R}^3\times Y\times Y} A' H\left(\frac{B}{A'}\right) f'\1_{\{y<y'\}}dydy'dx\\\\
& & - & \int\!\!\!\int\!\!\!\int_{\mathbb{R}^3\times Y\times Y} B f'\frac{d_u H(f')}{2}\1_{\{y<y'\}}dydy'dx.
\end{array}
\end{equation}
\begin{flushright}
$\Box$ 
\end{flushright}

\begin{Lemma} Forall $T>0$, there exists a constant $C_T>0$ depending only on $T$, the initial values $N_0$,$M_0$,$E_0$ and the truncature parameter $C_0$ such that for all $t\in [0,T]$, 
\begin{equation}\label{bornes Lp}
\begin{array}{lll}
\int\!\!\!\int_{\mathbb{R}^3\times Y}f(t,x,y)^s dxdy & \leq & \int\!\!\!\int_{\mathbb{R}^3\times Y}f^0(x,y)^s dxdy+C_T\\\\
& & -\frac{1}{2}\int_0^T\!\!\!\int\!\!\!\int\!\!\!\int_{\mathbb{R}^3\times Y\times Y} A\sup(f,f^{\star})\inf(f,f^{\star})^s dy^{\star}dydxd\tau\\\\
& & -\left|\frac{s-\delta}{2}\right|\int_0^T\!\!\!\int\!\!\!\int\!\!\!\int_{\mathbb{R}^3\times Y\times Y} B \left(f'\right)^s\1_{\{y<y'\}}dy'dydxd\tau,
\end{array}
\end{equation}

\noindent where $s$ and $\delta$ are given by (\ref{comparaison des noyaux}).
\end{Lemma}

\noindent\underline{Proof}: We use the previous lemma with $H(u)=u^s$. We obtain
$$\begin{array}{lllll}
\int\!\!\!\int_{\mathbb{R}^3\times Y} f(t,x,y)^s dxdy  & \leq & & \hspace{-5mm}\int\!\!\!\int_{\mathbb{R}^3\times Y} f^0(x,y)^s dxdy\\\\ 
& & - & \frac{1}{2}\int_0^t\!\!\!\int\!\!\!\int\!\!\!\int_{\mathbb{R}^3\times Y\times Y} \!\!\!\!\!\!\!A\sup(f,f^{\star})\inf(f,f^{\star})^s dy^{\star}dydxd\tau\\\\
&  & + & \int_0^t\!\!\!\int\!\!\!\int\!\!\!\int_{\mathbb{R}^3\times Y\times Y}\left(\frac{B^s}{A'^{s-1}}\right) f'\1_{\{y<y'\}}dydy'dxd\tau\\\\
&  & - & \frac{s}{2}\int_0^t\!\!\!\int\!\!\!\int_{\mathbb{R}^3\times Y} B_1(y') (f')^s dy'dxd\tau. 
\end{array}$$

\noindent According to (\ref{comparaison des noyaux}) and (\ref{premier r\'esum\'e des estimations}),

$$\begin{array}{llll}
& \hspace{-15mm}\int_0^t\!\!\!\int\!\!\!\int\!\!\!\int_{\mathbb{R}^3\times Y\times Y}\left(\frac{B^s}{A'^{s-1}}\right) f'\1_{\{y<y'\}}dydy'dxd\tau\\\\
\leq & \int_0^t\!\!\!\int\!\!\!\int_{\mathbb{R}^3\times Y}\left(1+m'+\frac{|p'|^2}{2m'}+e'\right) f'dy'dxd\tau+\frac{1}{2}\int_0^t\!\!\!\int\!\!\!\int_{\mathbb{R}^3\times Y}B_1(y')^{\delta} f'dy'dxd\tau\\\\
\leq & TK_T +\frac{1}{2}\int_0^t\!\!\!\int\!\!\!\int_{\mathbb{R}^3\times Y}B_1(y')^{\delta} f'^{s\delta}f'^{1-s\delta}dy'dxd\tau. 
\end{array}
$$

\noindent We apply the Young inequality again with the exponent $1/\delta>1$:
$$\begin{array}{llll}
B_1(y')^{\delta} f'^{s\delta}f'^{1-s\delta} & \leq & \frac{\left(B_1(y')^{\delta} f'^{s\delta}\right)^{1/\delta}}{1/\delta}+\frac{\left(f'^{1-s\delta}\right)^{\left(1/\delta\right)^{\star}}}{\left(1/\delta\right)^{\star}}\\\\
& = & \delta B_1(y')f'^s+(1-\delta)f'^{\frac{1-s\delta}{1-\delta}}.  
\end{array}
$$

\noindent Thus we deduce 
$$\begin{array}{lllll}
\int\!\!\!\int_{\mathbb{R}^3\times Y} f(t,x,y)^s dxdy  & \leq & & \hspace{-5mm}\int\!\!\!\int_{\mathbb{R}^3\times Y} f^0(x,y)^s dxdy\\\\ 
& & - & \frac{1}{2}\int_0^t\!\!\!\int\!\!\!\int\!\!\!\int_{\mathbb{R}^3\times Y\times Y} \!\!\!\!\!\!\!A\sup(f,f^{\star})\inf(f,f^{\star})^s dy^{\star}dydxd\tau\\\\
& & + & TK_T +\frac{1-\delta}{2}\int_0^t\!\!\!\int\!\!\!\int_{\mathbb{R}^3\times Y}f'^{\frac{1-s\delta}{1-\delta}}dy'dxd\tau\\\\ 
& & + & \frac{\delta-s}{2}\int_0^t\!\!\!\int\!\!\!\int_{\mathbb{R}^3\times Y} B_1(y') (f')^s dy'dxd\tau.\\\\
\end{array}$$

\noindent We can use (\ref{bornes L^q, q<1}) since $\frac{1-s\delta}{1-\delta}\in\,]\,5/6\,,\,1\,[.$
$$\begin{array}{lllll}
\int\!\!\!\int_{\mathbb{R}^3\times Y} f(t,x,y)^s dxdy &  \leq & & \hspace{-5mm}\int\!\!\!\int_{\mathbb{R}^3\times Y} f^0(x,y)^s dxdy+C_T\\\\
& & - & \frac{1}{2}\int_0^t\!\!\!\int\!\!\!\int\!\!\!\int_{\mathbb{R}^3\times Y\times Y} \!\!\!\!\!\!\!A\sup(f,f^{\star})\inf(f,f^{\star})^s dy^{\star}dydxd\tau\\\\
& & + & \frac{\delta-s}{2}\int_0^t\!\!\!\int\!\!\!\int_{\mathbb{R}^3\times Y} B_1(y') (f')^s dy'dxd\tau.\\\\ 
\end{array}$$

\noindent We conclude noting that $\delta<\frac{1}{6s-5}<1<s.$
\begin{flushright}
$\Box$ 
\end{flushright}

\end{subsection}

\end{section}

\begin{section}{A stability result}

\noindent The proof of theorem \ref{Theoreme d'existence} relies on a stability theorem, which claims that we can pass to the limit in the equation $(ECF)$ in a certain sense, namely in an integral formulation.  

\begin{Defin} Let $T>0$ and let $f^0$ be a nonnegative initial data which satisfies (\ref{bornes sur la donn\'ee initiale}). A weak solution of (\ref{probl\`eme de Cauchy}) is a nonnegative function $f \in C([0,T],L^1(\mathbb{R}^3\times Y))$, verifying the estimates (\ref{bornes L^1, L^s et moments}) and (\ref{autres quantit\'es born\'ees}), satisfying $(ECF)$ in $\mathcal{D}'(]0,+\infty[\times\mathbb{R}^3\times Y)$, and such that $f(0)=f^0$. 
\end{Defin}

\noindent Now let us state the result we will prove in this section:

\begin{Th}\label{theoreme de stabilite} Let $(f_n)_{n\geq 1}$ be a sequence of weak solutions of (\ref{probl\`eme de Cauchy}), with initial data $f^0_n$, and such that

\begin{equation}\label{r\'egularit\'e des f_n}
\forall n\in\mathbb{N},\qquad f_n\in W^{1,1}(]0,+\infty[\times\mathbb{R}^3\times Y),
\end{equation}

\begin{equation}\label{bornes L^1, L^s et moments unif}
\sup_{n\geq 1}\sup_{t\in[0,T]}\int\!\!\!\int_{\mathbb{R}^3\times Y}\left(\left(1+m+\frac{|p|^2}{2m}+e+m|x|^2\right)f_n(t,x,y)+f_n(t,x,y)^q\right)dxdy\leq K_T,
\end{equation}

\noindent for all the exponents $q\in\,]\,5/6\,,\,s\,],$ and also

\begin{equation}\label{autres quantit\'es born\'ees unif}
\sup_{n\geq 1}\int_0^T\!\!\!\int_{\mathbb{R}^3} \left(D_1(f_n(t,x))+D_2(f_n(t,x))\right)dxdt\leq K_T.
\end{equation}

\noindent (the \textit{a priori} estimates hold uniformly in $n$).   

\noindent Then, up to a subsequence, $f_n\stackrel{}{\rightharpoonup} f$ weakly in $L^1(]0,T[\times \mathbb{R}^3_{loc}\times Y)$, where $f$ is a renormalized solution of (ECF). Furthermore, $f\in C([0,T],L^1(\mathbb{R}^3\times Y))$.
\end{Th}
%\end{subsection}

\begin{subsection}{Weak compactness of $(f_n)$}

\noindent Let $0<T<\infty$. The bounds on $f_n$ provides some weak compactness, and thus the existence of a limit $f$ after extraction. 
  
\begin{Lemma}\label{compacit\'e faible des f_n}
For all $R>0$, the sequence $(f_n)_{n\geq 1}$ is weakly compact in $L^1(]0,T[\times B_R\times Y)$.
\end{Lemma}

\noindent\underline{Proof}: We set $\Phi(\xi):=\xi^s$ and $\Psi\left(m,p,e\right):=m+\frac{|p|^2}{2m}+e$.
The function $\Phi$ is nondecreasing, nonnegative and ${\Phi(\xi)}/{\xi}\underset{\xi\rightarrow +\infty}{\longrightarrow} +\infty$, $\Psi$ is nonnegative and $\Psi(y)\underset{|y|\rightarrow +\infty}{\longrightarrow} +\infty$. The estimate (\ref{bornes L^1, L^s et moments unif}) gives 
$$\sup_{n\geq 1}\int_0^T\!\!\!\int_{B_R}\!\int_{Y}\left(\left(1+\Psi(y)\right)f_n+\Phi\left(f_n\right)\right)dtdxdy< +\infty,$$
and we conclude by Dunford-Pettis theorem.
\begin{flushright}
$\Box$ 
\end{flushright}

\noindent Thus, there exists a nonnegative function $f$ such that for all $R>0$, $f_n\rightharpoonup f$ in \\${L^1(]0,T[\times B_R\times Y)}$ for a subsequence (not relabeled). Moreover, we can show easily (diagonal extraction) that the subsequence is not depending on $R$. Then we notice that in fact, $f\in L^{\infty}(]0,T[,L^1(\mathbb{R}^3\times Y))$ and 
\begin{equation}\label{estimation limite}
a.e\,\,t\in]0,T[,\qquad\int\!\!\!\int_{\mathbb{R}^3\times Y}\left(1+m+\frac{|p|^2}{2m}+e+m|x|^2\right)f(t,x,y)dxdy\leq K_T. 
\end{equation}

\noindent Moreover, since the function $\xi\mapsto|\xi|^s$ is convex, we have 
\begin{equation}\label{estimation limite 2}
a.e\,\,t\in]0,T[,\qquad\int\!\!\!\int_{\mathbb{R}^3\times Y}f(t,x,y)^sdxdy\leq K_T.  
\end{equation}

\end{subsection}

\begin{subsection}{Weak compactness of the renormalized coalescence term}\label{section compacite faible coag}

\noindent The bounds on $f_n$ are not enough to define the term $Q_c^-(f_n,f_n)$ as a distribution, unlike the term $Q_c^+(f_n,f_n)$. So, it seems that renormalization is necessary to obtain well-defined and weakly compact coalescence terms. 

\begin{Lemma}\label{compacit\'e faible de Q_c+ brut}
For all $R>0$, the sequence $(Q_c^+(f_n,f_n))_{n\geq 1}$ is weakly compact in\\ $L^1(]0,T[\times B_R\times Y_R)$, where$\vspace{1mm}$ $Y_R:=]0,R[\times B_R\times]0,R[$.
\end{Lemma}

\noindent\underline{Proof}: Let $E$ be a measurable subset of $]0,T[\times B_R\times Y_R$. We set $\varphi(t,x,y):=\1_{E}(t,x,y)$. Performing the change of variables $(y,y^{\star})\rightarrow(y^{\star},y-y^{\star})$, we obtain

$$\begin{array}{llll}
& \!\!\!\!\!\!\int_{Y_R}Q_c^+(f_n,f_n)(y)\varphi(y)dy = \frac{1}{2}\int_{Y_R}\!\int_{Y} A(y^{\star},y-y^{\star})f_n(y^{\star})f_n(y-y^{\star})\varphi(y)\1_{\{y^{\star}< y\}}dy^{\star}dy\\\\
= & \frac{1}{2}\underset{\begin{array}{c} \scriptstyle 0<m<R\\\scriptstyle p\in\mathbb{R}^3\\
\scriptstyle 0<e<R\end{array}}{\int\!\!\!\int\!\!\!\int}
\underset{\begin{array}{c}\scriptstyle 0<m^{\star}<R-m\\\scriptstyle p^{\star}\in B(-p,R)\\\scriptstyle 0<e^{\star}<R-e-E_-(m,m^{\star},p,p^{\star})\end{array}}{\int\!\!\!\int\!\!\!\int} 
A(y,y^{\star})f_n(y)f_n(y^{\star})\varphi(y+y^{\star})dy^{\star}dy.\\\\
\end{array}$$

\noindent In fact, we only integrate over $p\in B_{2R}$ because 
$$\frac{|p|^2}{2m}\leq \frac{|p|^2}{2m}+\frac{|p^{\star}|^2}{2m^{\star}}=\frac{|p+p^{\star}|^2}{2(m+m^{\star})}+E_-(m,m^{\star},p,p^{\star})\leq \frac{R^2}{2(m+m^{\star})}+R,$$
which yields
$$|p|^2\leq R^2\frac{m}{m+m^{\star}}+2mR\leq3R^2$$
(and the same applies to $p^{\star}$ because of the symmetry in the previous computation).

\noindent Thus we have 
$$\begin{array}{llll}
 \int_{Y_R}Q_c^+(f_n,f_n)(y)\varphi(y)dy 
= \frac{1}{2}\underset{\begin{array}{c} \scriptstyle 0<m<R\\\scriptstyle p\in B_{2R}\\
\scriptstyle 0<e<R\end{array}}{\int\!\!\!\int\!\!\!\int}
\underset{\begin{array}{c}\scriptstyle 0<m^{\star}<R-m\\\scriptstyle p^{\star}\in B(-p,R)\\\scriptstyle 0<e^{\star}<R-e-E_-(m,m^{\star},p,p^{\star})\end{array}}{\int\!\!\!\int\!\!\!\int}
Af_nf_n^{\star}\varphi'dy^{\star}dy.\\\\
\end{array}$$

\noindent Using the inequality 
\begin{equation}\label{decompostion inf sup}
Af_nf_n^{\star}\leq \frac{1}{M^{s-1}}A\sup(f_n,f_n^{\star})\inf(f_n,f_n^{\star})^s\1_{\{\inf(f_n,f_n^{\star})>M\}}+MA\sup(f_n,f_n^{\star})\1_{\{\inf(f_n,f_n^{\star})\leq M\}},
\end{equation}

\noindent we obtain

$$\begin{array}{lllll}
\int_{Y_R}\!\!\! Q_c^+(f_n,f_n)(y)\varphi(y)dy  
& \leq & & \frac{D_1(f_n(t,x))}{M^{s-1}}+\frac{M}{2}\int_{Y_{2R}}\int_{Y_{2R}}A\sup(f_n,f_n^{\star})\varphi'dy^{\star}dy\\\\
& \leq & & \frac{D_1(f_n(t,x))}{M^{s-1}}+M\int_{Y_{2R}}\int_{Y_{2R}}A f_n\varphi'dy^{\star}dy\\\\
& \leq & & \frac{D_1(f_n(t,x))}{M^{s-1}}+M \|A\|_{\infty, {Y^2_{2R}}}\int_{Y_{2R}}\int_{Y_{2R}} f_n\varphi'dy^{\star}dy.
\end{array}$$

\noindent We fix $\eps>0$ and choose $M$ such that $1/M^{s-1}\leq\eps$.
 
\noindent So, we can write
$$\begin{array}{llllll}
\int_{Y_R}\!\!\! Q_c^+(f_n,f_n)(y)\varphi(y)dy  
& \leq & & \!\!\!\!\!\!\!\eps D_1(f_n(t,x)) + M\|A\|_{\infty, {Y^2_{2R}}}\int_{Y_{2R}}\int_{Y_{2R}}f_n\varphi'dy^{\star}dy\\\\
& \leq & & \!\!\!\!\!\!\!\eps D_1(f_n(t,x))+M\|A\|_{\infty, {Y^2_{2R}}}\int_{Y}\int_{Y}f_n\varphi^{\star}dy^{\star}dy.
\end{array}$$

\noindent Eventually, in view of (\ref{bornes L^1, L^s et moments unif}) and (\ref{autres quantit\'es born\'ees unif}), we obtain
$$\begin{array}{lllll}
\int_0^T\!\!\!\int_{B_R}\int_{Y_R}\!\!\! Q_c^+(f_n,f_n)\varphi(t,x,y)dydxdt  
& \leq & & \hspace{-18mm} M\|A\|_{\infty, {Y^2_{2R}}}\!\!\int_0^T\!\!\!\int_{B_R}\!\int_{Y}\!\int_{Y}f_n\varphi^{\star}dy^{\star}dydxdt
\\\\ 
& & + \hspace{3mm}\eps K_T.
\end{array}$$

\noindent We conclude by letting $mes(E)\rightarrow 0$ and using the weak compactness of $(f_n).$
\begin{flushright}
$\Box$ 
\end{flushright}

\begin{Cor}\label{compacit\'e faible de Q_c+}
For all $R>0$, the sequence $\left(\frac{Q_c^+(f_n,f_n)}{1+f_n}\right)_{n\geq 1}$ is weakly compact in $L^1(]0,T[\times B_R\times Y_R)$, where$\vspace{1mm}$ $Y_R:=]0,R[\times B_R\times]0,R[$.
\end{Cor}

\noindent\underline{Proof}: It's obvious by the previous lemma and Dunford-Pettis theorem since\\ $\frac{Q_c^+(f_n,f_n)}{1+f_n}\leq Q_c^+(f_n,f_n).$
\begin{flushright}
$\Box$ 
\end{flushright}

\begin{Lemma}\label{convergence faible de L}
For all $R>0$, $Lf_n\rightharpoonup Lf$ weakly in $L^1(]0,T[\times B_R\times Y_R)$, where$\vspace{1mm}$ $Y_R:=]0,R[\times B_R\times]0,R[$.
\end{Lemma}

\noindent\underline{Proof}: Let $\varphi(t,x,y)\in L^{\infty}(]0,T[\times B_R\times Y_R).$ We have
$$\begin{array}{llll}
& \int_{Y_R}\!\!\! Lf_n(y)\varphi(y)dy = \int_{Y_R}\!\int_{Y} A(y,y^{\star})f_n(y^{\star})\varphi(y) dy^{\star}dy.
\end{array}$$

\noindent We fix $\eps>0$ and, in view of the assumptions (\ref{symetrie A}) and (\ref{asymptotique2 de A}), we choose $R^{\star}>0$ such that
$$\forall |y^{\star}|>R^{\star},\qquad \int_{Y_R}\frac{A(y,y^{\star})}{|y^{\star}|}dy\leq \eps.$$

\noindent We can write
$$\begin{array}{llll}
\int_0^T\int_{B_R}\int_{Y_R}\!\!\!Lf_n\,\varphi dydxdt & = & & \hspace{-7mm}\int_0^T\int_{B_R}\int_{Y_R}\!\int_{Y_{R^{\star}}}A(y,y^{\star})f_n(y^{\star})\varphi(y) dy^{\star} dydxdt\\\\
& & + &\int_0^T\int_{B_R}\int_{Y_R}\!\int_{Y-Y_{R^{\star}}}A(y,y^{\star})f_n(y^{\star})\varphi(y) dy^{\star} dydxdt.\\\\
\end{array}$$
\noindent First, 
$$\int_0^T\int_{B_R}\int_{Y_R}\!\int_{Y_{R^{\star}}}Af_n(y^{\star})\varphi(y) dy^{\star} dydxdt\underset{n}{\longrightarrow} \int_0^T\int_{B_R}\int_{Y_R}\!\int_{Y_{R^{\star}}}Af(y^{\star})\varphi(y) dy^{\star} dydxdt.$$
\noindent Indeed, setting $\theta(t,x,y^{\star})=\int_{Y_R}A(y,y^{\star})\varphi(t,x,y)dy$, we have
$$\int_0^T\int_{B_R}\int_{Y_R}\!\int_{Y_{R^{\star}}}Af_n(y^{\star})\varphi(y) dy^{\star} dydxdt=\int_0^T\int_{B_R}\int_{Y_{R^{\star}}}\theta(t,x,y^{\star})f_n(t,x,y^{\star})dy^{\star}dxdt$$  
\noindent and we conclude by lemma \ref{compacit\'e faible des f_n}, because the assumption (\ref{A loc born\'ee}) implies\\ $\theta\in L^{\infty}(]0,T[\times B_R\times Y_{R^{\star}})$.

\noindent Moreover,
$$\begin{array}{lll}
\hspace{-7mm}\left|\int_0^T\int_{B_R}\int_{Y_R}\!\int_{Y-Y_{R^{\star}}}Af_n(y^{\star})\varphi(y) dy^{\star} dydxdt\right|
 \leq \eps\,\|\varphi\|_{\infty}\int_0^T \int_{B_R}\!\int_{Y-Y_{R^{\star}}}\!\!\!\!|y^{\star}|f_n(y^{\star})dy^{\star}dxdt\\\\
\hphantom{\hspace{-7mm}\left|\int_0^T\int_{B_R}\int_{Y_R}\!\int_{Y-Y_{R^{\star}}}Af_n(y^{\star})\varphi(y) dy^{\star} dydxdt\right|
 }\leq  \eps T\,\|\varphi\|_{\infty}\,K_T,
\end{array}$$

\noindent and the inequality (\ref{estimation limite}) yields

$$\begin{array}{lll}
\hspace{-7mm}\left|\int_0^T\int_{B_R}\int_{Y_R}\!\int_{Y-Y_{R^{\star}}}Af(y^{\star})\varphi(y) dy^{\star} dydxdt\right|
\leq \,\,\eps\,T\,\|\varphi\|_{\infty}\,K_T.
\end{array}$$

\noindent Finally, we infer
$$\begin{array}{llll}
\left|\int_0^T\int_{B_R}\int_{Y_R}\!\!\!Lf_n\,\varphi dydxdt - \int_0^T\int_{B_R}\int_{Y_R}\!\!\!Lf\,\varphi dydxdt\right|\leq \underset{n\rightarrow +\infty}{o(1)}+C(T,R,\varphi)\eps.
\end{array}$$ 

\begin{flushright}
$\Box$ 
\end{flushright}

\begin{Cor}\label{compacit\'e faible de Q_c-}
For all $R>0$, the sequence $\left(\frac{Q_c^-(f_n,f_n)}{1+f_n}\right)_{n\geq 1}$ is weakly compact in $L^1(]0,T[\times B_R\times Y_R)$, where$\vspace{1mm}$ $Y_R:=]0,R[\times B_R\times]0,R[$.
\end{Cor}

\noindent\underline{Proof}: It's obvious because $\frac{Q_c^-(f_n,f_n)}{1+f_n}=\frac{f_n}{1+f_n}Lf_n\leq Lf_n.$
\begin{flushright}
$\Box$ 
\end{flushright}

\end{subsection}

\begin{subsection}{Weak convergence of the fragmentation term}\label{section compacite faible frag}

\noindent Since they are linear, the fragmentation terms easily pass to the limit, and we have the following lemma.

\begin{Lemma}\label{convergence faible fragmentation}
For all $R>0$, we have 
\begin{itemize}
 \item[(i)] $Q_f^+(f_n)\rightharpoonup Q_f^+(f)$ weakly in $L^1(]0,T[\times B_R\times Y_R),$
\item[(ii)] $Q_f^-(f_n)\rightharpoonup Q_f^-(f)$ weakly in $L^1(]0,T[\times B_R\times Y_R),$
\end{itemize}
\noindent where$\vspace{1mm}$ $Y_R:=]0,R[\times B_R\times]0,R[$.
\end{Lemma}

\noindent\underline{Proof}: The part (ii) results immediately from (\ref{B_1 loc born\'ee}), and the proof of (i) is the same as lemma \ref{convergence faible de L}.
\begin{flushright}
$\Box$ 
\end{flushright}
\end{subsection} 

\begin{subsection}{Strong compactness of $y$-averages}

\noindent Strong compactness is needed to pass to the limit in coalescence terms (because they are quadratic), that's why we use the following averaging lemma, inspired by \cite{mischler1}, \cite{diperna-lions}, and \cite{diperna-lions-meyer}: 

\begin{Th}\label{lemme de moyenne} Let $(g_n)$ be a bounded sequence in $L^1(]0,T[\times \mathbb{R}^3\times Y)$ and weakly compact in $L^1(]0,T[\times B_R\times Y_R)$, for all $R>0$. Let $(G_n)$ be a bounded sequence in $L^1(]0,T[\times B_R\times Y_R)$ for all $R>0$. We assume that 
$$\partial_{t} g_n+\frac{p}{m}.\nabla_x g_n =G_n\quad\mbox{in}\quad \mathcal{D}'(]0,+\infty[\times\mathbb{R}^3\times Y).$$

\noindent Then, for any function $\Psi\in L^{\infty}(Y^2)$, with compact support, the sequence 
$$\left(\int_Y g_n(t,x,y)\Psi(y,y^{\star})dy\right)_{n\in\mathbb{N}}$$ 
\noindent is strongly compact in $L^1(]0,T[\times B_R\times Y_R)$, for all $R>0$.
\end{Th}\vspace{2mm}

\noindent This result can be improved:
\begin{Cor}\label{corollaire 1 du lemme de moyenne} With the assumptions of theorem \ref{lemme de moyenne}, we also have: 
 
\noindent for all $R>0$ and for any function $\Psi\in L^{\infty}(]0,T[\times B_R\times Y_R^2)$, the sequence 
$$\left(\int_Y g_n(t,x,y)\Psi(t,x,y,y^{\star})dy\right)_{n\in\mathbb{N}}$$ 
\noindent is strongly compact in $L^1(]0,T[\times B_R\times Y_R)$.
\end{Cor}

\noindent\underline{Proof}: The case of separated variables is obvious. Then, we proceed by a density argument as in \cite{diperna-lions}.
\begin{flushright}
$\Box$ 
\end{flushright}

\begin{Cor}\label{corollaire 2 du lemme de moyenne} With the assumptions of theorem \ref{lemme de moyenne}, we also have: for all $R>0$ and for any sequence $(\Psi_n)$ bounded in $L^{\infty}(]0,T[\times B_R\times Y_R)$ which converges a.e to $\Psi\in L^{\infty}(]0,T[\times B_R\times Y_R),$ the sequence 
$$\left(\int_Y g_n(t,x,y)\Psi_n(t,x,y)dy\right)_{n\in\mathbb{N}}$$ 
\noindent is strongly compact in $L^1(]0,T[\times B_R)$.
\end{Cor}

\noindent\underline{Proof}: Let $\eps>0$. The sequence $(g_n)$ being weakly compact in $L^1(]0,T[\times B_R\times Y_R)$, there exists $\delta>0$ such that 
$$\forall E\in \mathcal{B}(]0,T[\times B_R\times Y_R),\quad |E|<\delta,\qquad\sup_n \int\!\!\!\int\!\!\!\int_E |g_n|dtdxdy\leq \eps.$$

\noindent Then, by Egoroff theorem, there exists $E_0\in\mathcal{B}(]0,T[\times B_R\times Y_R)$ such that $|E_0|<\delta$ and $\Psi_n$ converge uniformly to $\Psi$ on $E_1:=(]0,T[\times B_R\times Y_R)\setminus E_0.$ Whence

$$\begin{array}{lllll}
\left\|\int_{Y_R} g_n \Psi_n dy-\int_{Y_R} g_n \Psi dy\right\|_{L^1(]0,T[\times B_R)} & \leq & \int_0^T\!\!\!\int_{B_R}\!\int_{Y_R}\!|g_n||\Psi_n-\Psi|dydxdt\\\\
 & \leq  & 2 C\,\,\eps +\sup_{E_1}|\Psi_n-\Psi|\int\!\!\!\int\!\!\!\int_{E_1}|g_n|dydxdt\\\\
& = & 2C\,\,\eps+\underset{n\rightarrow+\infty}{o(1)}.
  \end{array}
$$

\noindent We infer $$\left\|\int_{Y_R} g_n \Psi_n dy-\int_{Y_R} g_n \Psi dy\right\|_{L^1(]0,T[\times B_R)}\overset{n}{\longrightarrow} 0.$$

\noindent The sequence $\left(\int_{Y_R}g_n\Psi dy\right)$ being compact in $L^1(]0,T[\times B_R)$ in view of corollary \ref{corollaire 1 du lemme de moyenne}, the results follows.
\begin{flushright}
$\Box$ 
\end{flushright}

\noindent Now, we are able to establish the strong compactness of the sequence of $f_n$ $y$-averages, and also the $(Lf_n)$ one. 

\begin{Lemma}\label{convergence L^1 des moyennes des f_n} For all $R>0$, and for all function $\Psi\in L^{\infty}(Y)$ with compact support,
$$\int_Y f_n(t,x,y)\Psi(y)dy\overset{n}{\longrightarrow} \int_Y f(t,x,y)\Psi(y)dy\quad \mbox{in}\quad L^1(]0,T[\times B_R).$$ 
\end{Lemma}

\noindent\underline{Proof}: Since it is not clear that $\left(Q_c^-(f_n,f_n)\right)_n$ is bounded in $L^1(]0,T[\times B_R\times Y_R)$, we can not directly apply theorem \ref{lemme de moyenne} to the sequence $(f_n)$. 
For $\nu>0$, we consider the sequence $g^{\nu}_n:=\frac{1}{\nu}\log(1+\nu f_n)$ and we set 
$$G_n^{\nu}:=\frac{Q_c^+(f_n,f_n)}{1+\nu f_n}-\frac{Q_c^-(f_n,f_n)}{1+\nu f_n}+\frac{Q_f^+(f_n)}{1+\nu f_n}-\frac{Q_f^-(f_n)}{1+\nu f_n}.$$ 

\noindent By the assumptions on $(f_n)$, we have 
\begin{equation}\label{edp sur les g_n^delta}
 \partial_{t} g_n^{\nu}+\frac{p}{m}.\nabla_x g_n^{\nu} =G_n^{\nu}\quad\mbox{in}\quad \mathcal{D}'(]0,+\infty[\times\mathbb{R}^3\times Y).
\end{equation}

\noindent Since $0\leq g_n^{\nu} \leq f_n$, the weak compactness of $(f_n)$ established in the lemma \ref{compacit\'e faible des f_n} implies that $(g_n^{\nu})$ is also weakly compact. Similarly, the sequence $(g_n^{\nu})$ is bounded in $L^1(]0,T[\times \mathbb{R}^3\times Y)$. 
Then, by corollaries \ref{compacit\'e faible de Q_c+}, \ref{compacit\'e faible de Q_c-} and lemma \ref{convergence faible fragmentation}, $(G_n^{\nu})$ is bounded in $L^1(]0,T[\times B_R\times Y_R)$. Therefore, theorem \ref{lemme de moyenne} applies to $(g_n^{\nu})$ for all $\nu>0$. In particular, for all function $\Psi\in L^{\infty}(Y)$ with compact support and for all $\nu>0$, the sequence\\
$\left(\int_Y g_n^{\nu}(t,x,y)\Psi(y)dy\right)_{n}$ is compact in $L^1(]0,T[\times B_R),$ thus, by the uniqueness of weak limit,

\begin{equation}\label{convergence L^1 des moyennes des g_n^delta}
\int_Y g_n^{\nu}(t,x,y)\Psi(y)dy\overset{n}{\longrightarrow} \int_Y g^{\nu}(t,x,y)\Psi(y)dy\quad \mbox{in}\quad L^1(]0,T[\times B_R), 
\end{equation}

\noindent where $g^{\nu}$ is the weak limit of $(g_n^{\nu})$ (up to an extraction).

\noindent The result follows because 
\begin{equation}\label{convergence unif en n selon delta}
\sup_n \sup_{t\in [0,T]}\int\!\!\!\int_{\mathbb{R}^3\times Y} |g_n^{\nu}-f_n|dydx \underset{\nu\rightarrow 0}{\longrightarrow} 0,
\end{equation}

\noindent which implies the strong compactness in $L^1(]0,T[\times B_R)$ of the sequence\\ $\left(\int_Y f_n(t,x,y)\Psi(y)dy\right)_n.$

\noindent To show (\ref{convergence unif en n selon delta}), we can use the inequality
\begin{equation}\label{Taylor log}
\forall M>0,\qquad0\leq u-\frac{1}{\nu}\log(1+\nu u)=\frac{\nu M}{2}u\1_{\{u\leq M\}}+u\1_{\{u>M\}}.
\end{equation}
 
\noindent Then we obtain, for all $n$ and for all $t\in[0,T]$, 
$$\begin{array}{lll}
\int\!\!\!\int_{\mathbb{R}^3\times Y} |g_n^{\nu}-f_n|dydx & \leq & \frac{\nu M}{2}\int\!\!\!\int_{\mathbb{R}^3\times Y}f_ndydx+\int\!\!\!\int_{\mathbb{R}^3\times Y}f_n\1_{\{f_n>M\}}dydx\\\\
& \leq & \frac{\nu M}{2}K_T \,+\, \frac{1}{M^{s-1}}\int\!\!\!\int_{\mathbb{R}^3\times Y}f_n^sdydx\leq \left(\frac{\nu M}{2}+\frac{1}{M^{s-1}}\right)K_T.
\end{array}$$

\noindent We conclude by letting $\nu\rightarrow0$, and $M\rightarrow+\infty.$ 
\begin{flushright}
$\Box$ 
\end{flushright}

\begin{Prop}We set $\rho_n(t,x):=\int_{Y}f_n(t,x,y)dy$ and $\rho(t,x):=\int_{Y}f(t,x,y)dy.$\\ 
\noindent Then, up to a subsequence, we have, for all $R>0$,  
\begin{equation}\label{convergence L^1 des moyennes}
\rho_n\longrightarrow\rho\quad \mbox{in}\quad L^1(]0,T[\times B_R)\quad\mbox{and}\quad a.e.
\end{equation}
\end{Prop}

\noindent\underline{Proof}: We have
$$\rho_n=\rho_n^M+\sigma_n^M,\quad\mbox{where}\quad \rho_n^M:=\int_{Y_M}f_n(t,x,y)dy.$$

\noindent By the preceding lemma, $\rho_n^M\underset{n}{\longrightarrow}\rho^M:=\int_{Y_M}f(t,x,y)dy$ in $L^1(]0,T[\times B_R)$ for all $M>0$,
and
$$\begin{array}{lll}
\sigma_n^M:=\int_{Y-Y_M}f_n(t,x,y)dy & \leq & \frac{1}{M}\int_{Y-Y_M}|y|f_n(t,x,y)dy\\\\
& \leq & \frac{Cte}{M}\int_Y \left(m+\frac{|p|^2}{2m}+e\right)f_n(t,x,y)dy,
\end{array}
$$

\noindent whence $\sigma_n^M\underset{M\rightarrow+\infty}{\longrightarrow}0$ in $L^1(]0,T[\times B_R)$, uniformly in $n$.
\begin{flushright}
$\Box$ 
\end{flushright}

\begin{Lemma}\label{Convergence forte de Lf_n} For all $R>0$, we have, up to a subsequence,
\begin{equation}\label{conv forte Lf_n}
Lf_n\longrightarrow Lf\quad\mbox{in}\quad L^1(]0,T[\times B_R\times Y_R)\quad\mbox{and}\quad a.e. 
\end{equation}
\end{Lemma}

\noindent\underline{Proof}: Applying the corollary \ref{corollaire 1 du lemme de moyenne} with $\Psi(y,y^{\star})=A(y,y^{\star})\1_{\{y\in Y_R\}}\1_{\{y^{\star}\in Y_{R^{\star}}\}}$, we infer that the sequence $\left(\int_{Y_{R^{\star}}} g_n^{\nu}(t,x,y^{\star})A(y,y^{\star})dy^{\star}\right)_n$ is compact in $L^1(]0,T[\times B_R\times Y_R)$ for all $R^{\star}>0$. Using (\ref{convergence unif en n selon delta}) again, we obtain, for all $R^{\star}>0,$ the compactness of $\left(\int_{Y_{R^{\star}}} f_n(t,x,y^{\star})A(y,y^{\star})dy^{\star}\right)_n$ in $L^1(]0,T[\times B_R\times Y_R)$. We conclude similarly as for the proof of lemma \ref{convergence faible de L}, establishing
$$\lim_{R^{\star}\rightarrow +\infty}\sup_{n}\left\|\int_{Y_{R^{\star}}}Af_n^{\star}dy^{\star}-\int_{Y}Af_n^{\star}dy^{\star}\right\|_{L^1(]0,T[\times B_R\times Y_R)}=0,$$

\noindent and identifying the weak limits.
\begin{flushright}
$\Box$ 
\end{flushright}

\end{subsection}

\begin{subsection}{Regularity in time of the limit $f$}

\noindent In this subsection, we show the continuity in time of the limit $f$, which gives a sense to the Cauchy data $f(0)=f^0$.

\begin{Prop} In fact, we have $f\in C([0,T],L^1(\mathbb{R}^3\times Y))$. 
\end{Prop}

\noindent\underline{Proof}: We use the integral formulation. Each $g_n^{\nu}$ is a distributional solution of the renormalized equation, by (\ref{edp sur les g_n^delta}), so a mild solution. Therefore we have, for $a.e\,\,(x,y)\in\mathbb{R}^3\times Y,$ and for all $t,t+h\in [0,T],$
$$g_n^{\nu\sharp}(t+h,x,y)-g_n^{\nu\sharp}(t,x,y)=\int_t^{t+h} G_n^{\nu}(\sigma,x,y)d\sigma,$$

\noindent thus

$$\|g_n^{\nu\sharp}(t+h)-g_n^{\nu\sharp}(t)\|_{L^1(B_R\times Y_R)}\leq\int\!\!\!\int_{B_R\times Y_R}\int_t^{t+h} \left|G_n^{\nu}(\sigma,x,y)\right|d\sigma.$$

\noindent Moreover, by the subsections \ref{section compacite faible coag} and \ref{section compacite faible frag}, the sequence $\left(G_n^{\nu}\right)_{n}$ is weakly compact in $L^1(]0,T[\times B_R\times Y_R),$ thus for all $t\in [0,T]$,
$$\lim_{h\rightarrow 0}\,\,\sup_n \|g_n^{\nu\sharp}(t+h)-g_n^{\nu\sharp}(t)\|_{L^1(B_R\times Y_R)}= 0.$$

\noindent Therefore the sequence $(g_n^{\nu\sharp})$ is equicontinuous in $C([0,T],L^1(B_R\times Y_R))$. By the compactness of $[0,T]$, this sequence is in fact uniformly equicontinuous, thus
$$\lim_{h\rightarrow 0}\,\,\sup_n \sup_{t\in [0,T]}\|g_n^{\nu\sharp}(t+h)-g_n^{\nu\sharp}(t)\|_{L^1(B_R\times Y_R)}= 0.$$

\noindent Then, (\ref{convergence unif en n selon delta}) and the estimate (\ref{bornes L^1, L^s et moments unif}) yield
$$\lim_{h\rightarrow 0}\,\,\sup_n \sup_{t\in [0,T]}\|f_n^{\sharp}(t+h)-f_n^{\sharp}(t)\|_{L^1(\mathbb{R}^3\times Y)}= 0.$$

\noindent Ascoli theorem entails that the sequence $(f_n^{\sharp})$ is compact in $C([0,T],L^1(\mathbb{R}^3\times Y))$.\\ The uniqueness of the limit in $\mathcal{D'}(]0,+\infty[\times\mathbb{R}^3\times Y)$ yields $f^{\sharp}\in C([0,T],L^1(\mathbb{R}^3\times Y)),$ and so $f\in C([0,T],L^1(\mathbb{R}^3\times Y))$ by change of variables.
\begin{flushright}
$\Box$ 
\end{flushright}

\end{subsection}

\begin{subsection}{Passing to the limit in a new integral equation}

\noindent Even if the renormalization provides weak compactness, a new problem appears: we will not be able to pass to the weak limit in $(ECFR)$, because of the non-linearity of the factor $f_n/(1+f_n)$. That's why we need another formulation to our problem, which avoids the renormalization. But, remember that the term $Q_c^-(f,f)$ can not be defined as a distribution, so we will use an integral equation which doesn't involve this term. We proceed as in \cite{gerard}.\vspace{5mm}

\noindent We denote by $T$ the linear transport operator
 $$T=\partial_t+\frac{p}{m}.\nabla_x.$$

\noindent Let $T^{-1}$ be the resolvant of transport operator, defined by: for $g(t,x,y),$ we set 
$u=T^{-1}g$ if $u|_{t=0}=0$ and $Tu=g.$ So, we have 
$$T^{-1}g(t,x,m,p,e):=\int_0^t g(s,x-(t-s)p/m,m,p,e)ds.$$ 

\noindent $T^{-1}$ satisfies the following properties:\vspace{2mm}

(i) Forall $R>0$, $T^{-1}\big(L^1(]0,T[\times B_R\times Y_{R})\big)\subset C([0,T],L^1(B_R\times Y_{R}))$ continuously and weakly continuously.

(ii) $T^{-1}$ is nonnegative ($\forall g\geq 0, T^{-1}g\geq 0$).\vspace{2mm}

\noindent Forall $F\in C([0,T], L^1(B_R\times Y_{R}))$ such that $TF\geq 0,$ we set $$T_F^{-1}=e^{-F}T^{-1}e^{F}.$$
This operator is well defined from $L^1(]0,T[\times B_R\times Y_{R})$ to $C([0,T], L^1(B_R\times Y_{R}))$ and has the same continuity properties as $T^{-1}$.

\noindent Moreover, if $(F_n)$ is a bounded sequence in $C([0,T], L^1(B_R\times Y_{R}))$ such that $TF_n\geq 0$, if $F_n(t,x,y)\rightarrow F(t,x,y)$ for all $t$ and $a.e\,(x,y)$, and if $g_n\rightharpoonup g$ weakly in $L^1(]0,T[\times B_R\times Y_R)$, then 
$$\forall t\in[0,T],\qquad T_{F_n}^{-1}g_n(t)\rightharpoonup T_{F}^{-1}g(t)\quad\mbox{weakly in}\quad L^1(B_R\times Y_R).$$

\noindent The operator $T_F^{-1}$ allows us to build an new formulation of our problem, which is better because it only involves $Q_c^+(f_n,f_n)$, $Q_f^+(f_n)$, $Q_c^-(f_n)$:

\begin{Lemma} $f\in C([0,T],L^1(\mathbb{R}^3\times Y))$ is a mild solution of $(ECF)$ with initial data $f(0)=f^0$ if and only if
\begin{equation}\label{formulation int\'egrale}
f=e^{-F}f^0(x-tp/m,y)+T_{F}^{-1}(Q_c^+(f,f))+T_{F}^{-1}(Q_f^+(f))-T_{F}^{-1}(Q_c^-(f)),
\end{equation}

\noindent where $F:=T^{-1}(Lf).$
\end{Lemma}

\noindent\underline{Proof}: The result is deduced from the following fact: if $f$ is a distributional solution of $(ECF)$, then
$$T(e^{F}f)=TF e^F f+e^F Tf=e^{F}(f\,Lf+Q_c^+(f,f)-Q_c^-(f,f)+Q_f^+(f)-Q_f^-(f))$$
$$\hphantom{T(e^{F}f)=TF e^F f+e^F Tf}=e^{F}(Q_c^+(f,f)+Q_f^+(f)-Q_f^-(f)).$$
\begin{flushright}
$\Box$ 
\end{flushright}

\noindent Now, we can finish the proof of theorem \ref{theoreme de stabilite}.\vspace{2mm}  

\noindent\underline {End of the proof of theorem \ref{theoreme de stabilite}}: We will pass to the weak limit in the following equation, satisfied by each $(f_n)$:

\begin{equation}\label{formulation int\'egrale pour les f_n}
f_n=e^{-F_n}f_n^0(x-tp/m,y)+T_{F_n}^{-1}(Q_c^+(f_n,f_n))+T_{F_n}^{-1}(Q_f^+(f_n))-T_{F_n}^{-1}(Q_c^-(f_n)),
\end{equation}

\noindent where $F_n:=T^{-1}(Lf_n).$ 

\noindent Notice that in view of (\ref{conv forte Lf_n}) and the continuity properties of $T^{-1}$, the sequence $(F_n)$ is bounded in $C([0,T], L^1(B_R\times Y_{R}))$, and $F_n(t,x,y)\rightarrow F(t,x,y)$ for all $t$ and $a.e\,(x,y)$. Thus, we can pass to the weak limit in the terms $T_{F_n}^{-1}(Q_f^+(f_n))$ and $T_{F_n}^{-1}(Q_c^-(f_n))$ thanks to the lemma \ref{convergence faible fragmentation}. The term $e^{-F_n}f_n^0(x-tp/m,y)$ can be treated with the continuity in $t=0$ (for the $L^1$-norm) of each $f_n$ and $f$, established in the previous section. Eventually, the last term $T_{F_n}^{-1}(Q_c^+(f_n,f_n))$ also pass to the weak limit thanks to the following lemma and proposition, which use the $a.e$ convergence of the $y$-averages obtained in the previous subsection.

\begin{Lemma}\label{Convergence intermediaire} For all $R>0$ and for all function $\varphi\in L^{\infty}(]0,T[\times B_R\times Y_R)$, we have, up to a subsequence,

$$\frac{\int_Y Q_c^+(f_n,f_n)(t,x,y)\varphi(t,x,y)dy}{1+\rho_n(t,x)}\overset{n}{\longrightarrow}\frac{\int_Y Q_c^+(f,f)(t,x,y)\varphi(t,x,y)dy}{1+\rho(t,x)}$$

\noindent in $L^1(]0,T[\times B_R)$ and a.e.
\end{Lemma}

\noindent\underline{Proof}: We have 

$$\frac{\int_Y Q_c^+(f_n,f_n)\varphi dy}{1+\rho_n}=\!\frac{1}{2}\int_Y f_n(t,x,y^{\star})\left(\!\int_Y  \frac{f_n(t,x,y)A(y,y^{\star})\varphi(t,x,y+y^{\star})dy}{1+\rho_n(t,x)}\right)dy^{\star}.$$

\noindent Now, we apply the corollary \ref{corollaire 1 du lemme de moyenne} with $\Psi(t,x,y,y^{\star})=A(y,y^{\star})\varphi(t,x,y+y^{\star})$ \\(notice that $\Psi\in L^{\infty}(]0,T[\times B_R\times Y_{2R}^2)$ thanks to (\ref{limitations})).\vspace{2mm} 

\noindent Therefore, the sequence $\left(\int_Y g_n^{\nu}(t,x,y)A(y,y^{\star})\varphi(t,x,y') dy\right)_n$ is compact in\\ $L^1(]0,T[\times B_R\times Y_R)$, and we have

\begin{equation}\label{convergence L^1 des moyennes des A g_n^delta}
\int_Y g_n^{\nu}(t,x,y)A(y,y^{\star})\varphi(t,x,y') dy\overset{n}{\longrightarrow} \int_Y g^{\nu}(t,x,y)A(y,y^{\star})\varphi(t,x,y') dy
\end{equation}

\noindent in $L^1(]0,T[\times B_R\times Y_R)$, for all $\nu>0$.\vspace{2mm}

\noindent Using $(\ref{convergence unif en n selon delta})$ again, we obtain, up to a subsequence, that

\begin{equation}\label{convergence L^1 des moyennes en A f_n}
\int_Y f_n(t,x,y)A(y,y^{\star})\varphi(t,x,y')dy\overset{n}{\longrightarrow} \int_Y f(t,x,y)A(y,y^{\star})\varphi(t,x,y')dy
\end{equation}

\noindent in $L^1(]0,T[\times B_R\times Y_R)$ and a.e.\vspace{2mm}

\noindent Up to another extraction, we infer, by (\ref{convergence L^1 des moyennes en A f_n}) and (\ref{convergence L^1 des moyennes}),
\begin{equation}\label{convergence pp intermediaire}
\int_Y \frac{f_n(t,x,y)A(y,y^{\star})\varphi(t,x,y')dy}{1+\rho_n(t,x)}\overset{n}{\longrightarrow}\int_Y\frac{ f(t,x,y)A(y,y^{\star})\varphi(t,x,y')dy}{1+\rho(t,x)}
\end{equation}

\noindent a.e in $(t,x,y)\in]0,T[\times B_R\times Y_R$.\vspace{2mm}

\noindent Applying the corollary \ref{corollaire 2 du lemme de moyenne} with $\Psi_n(t,x,y^{\star}):=\int_Y  \frac{f_n(t,x,y)A(y,y^{\star})\varphi(t,x,y')dy}{1+\rho_n(t,x)}$\\ (which satisfies the required assumptions because $\varphi$ is compact supported and $A$ is locally bounded), we obtain the compactness of the sequence\\ 
$\left(\int_Y g_n^{\nu}(t,x,y^{\star})\Psi_n(t,x,y^{\star})dy^{\star}\right)_{n\in\mathbb{N}}$ in $L^1(]0,T[\times B_R)$, and so, by (\ref{convergence unif en n selon delta}), we deduce that 
$\left(\int_Y f_n(t,x,y^{\star})\Psi_n(t,x,y^{\star})dy^{\star}\right)_{n\in\mathbb{N}}$ is compact.

\noindent Finally, we conclude that
$$\int_Y f_n(t,x,y^{\star})\Psi_n(t,x,y^{\star})dy^{\star}\overset{n}{\longrightarrow}\int_Y f(t,x,y^{\star})\Psi(t,x,y^{\star})dy^{\star}\quad\mbox{in}\quad L^1(]0,T[\times B_R).$$

\begin{flushright}
$\Box$ 
\end{flushright}

\begin{Prop} Up to a subsequence, we have, for all $R>0$, 
$$Q_c^+(f_n,f_n)\rightharpoonup Q_c^+(f,f)\quad\mbox{weakly in}\quad L^1(]0,T[\times B_R\times Y_R).$$
\end{Prop}

\noindent\underline{Proof}: We know by the lemma \ref{compacit\'e faible de Q_c+ brut} that there exists $\overline{Q}(t,x,y)$ such that for all $R>0$, $$Q_c^+(f_n,f_n)\rightharpoonup \overline{Q}\quad\mbox{weakly in}\quad L^1(]0,T[\times B_R\times Y_R).$$

\noindent By (\ref{convergence L^1 des moyennes}) and a standard integration argument (we can refer to \cite{mischler2} for a proof), it leads to 
$$\frac{Q_c^+(f_n,f_n)}{1+\rho_n}\rightharpoonup  \frac{\overline{Q}}{1+\rho}\quad\mbox{weakly in}\quad L^1(]0,T[\times B_R\times Y_R).$$

\noindent Moreover, the previous lemma shows that
$$\frac{Q_c^+(f_n,f_n)}{1+\rho_n}\rightharpoonup \frac{Q_c^+(f,f)}{1+\rho}\quad\mbox{weakly in}\quad L^1(]0,T[\times B_R\times Y_R).$$ 

\noindent We conclude identifying weak limits.
\begin{flushright}
$\Box$ 
\end{flushright}

\noindent We have shown that $f$ is a mild solution of $(ECF)$. Since $Q_c^+(f_n,f_n)$, $Q_f^+(f_n)$ and $Q_f^-(f_n)$ converge weakly to $Q_c^+(f,f)$, $Q_f^+(f)$ and $Q_f^-(f)$ respectively, these three terms lie in $L^1_{loc}$, and \textit{a fortiori}, 
$$\frac{Q_c^+(f,f)}{1+f}, \frac{Q_f^+(f)}{1+f}, \frac{Q_f^-(f)}{1+f}\in L^1_{loc}(]0,+\infty[\times\mathbb{R}^3\times Y).$$

\noindent The term $\frac{Q_c^-(f,f)}{1+f}$ is automatically in $L^1_{loc}$ because $Lf\in L^1(]0,T[\times B_R\times Y_R)$ for all $R>0$ and $\frac{Q_c^-(f,f)}{1+f}\leq Lf$.\\ Thus, $f$ is indeed a renormalized solution of $(ECF).$
\begin{flushright}
$\Box$ 
\end{flushright}

\end{subsection}
\end{section}
\vspace{10mm}

%\begin{center}
%\textbf{Acknowledgements} 
%\end{center}

%I would thank very much my advisors Pierre-Emmanuel Jabin and Florent Berthelin for their patience, their useful advices and especially for their benevolence.
%\clearpage

\end{document}